\documentclass[11pt]{article}

\usepackage{graphicx}
\usepackage[top=1.25in,bottom=1.25in,left=1in,right=1in]{geometry}
\usepackage{amsmath,amsfonts,amscd,amssymb,bm,epsfig,epsf,bbm,amssymb, prettyref} 
\usepackage{amsthm}
\usepackage{algorithm, algorithmic}
\usepackage[usenames]{color}
\usepackage{mathtools}

\usepackage{tikz}
\usepackage{pgfplots}
\usepackage{dsfont}
\usepackage[toc,page]{appendix}
\usepackage{makecell}
\usepackage{caption} 
\usepackage{subcaption}

\usepackage{natbib} 

\definecolor{darkred}{RGB}{100,0,0}
\definecolor{darkgreen}{RGB}{0,100,0}
\definecolor{darkblue}{RGB}{0,0,150}
\definecolor{red}{RGB}{255,0,0}

\usepackage{hyperref}
\hypersetup{colorlinks=true, linkcolor=darkred, citecolor=darkgreen, urlcolor=darkblue}
\usepackage{url}

\newtheorem{theorem}{Theorem}[section]
\newtheorem{lemma}[theorem]{Lemma}
\newtheorem{corollary}[theorem]{Corollary}

\newtheorem{definition}{Definition}

\newcommand{\E}{\operatorname{\mathbb{E}}}

\newcommand{\vct}[1]{\bm{#1}}
\newcommand{\mtx}[1]{\bm{#1}}

\newcommand{\colspan}{\operatorname{colspan}}

\newcommand{\rank}{\operatorname{rank}}

\numberwithin{equation}{section}

\definecolor{xli}{RGB}{200,50,50}

\numberwithin{equation}{section}
\DeclareMathOperator*{\argmin}{\arg\!\min}

\begin{document}

\title{Model-free Nonconvex Matrix Completion: Local Minima Analysis and Applications in Memory-efficient Kernel PCA}

\author{Ji Chen\thanks{Department of Mathematics, University of California Davis, Davis CA 95616}~~ and Xiaodong Li\thanks{Department of Statistics, University of California Davis, Davis CA 95616}}

\date{}
\maketitle

\begin{abstract}
This work studies low-rank approximation of a positive semidefinite matrix from partial entries via nonconvex optimization. We characterized how well local-minimum based low-rank factorization approximates a fixed positive semidefinite matrix without any assumptions on the rank-matching, the condition number or eigenspace incoherence parameter. Furthermore, under certain assumptions on rank-matching and well-boundedness of condition numbers and eigenspace incoherence parameters, a corollary of our main theorem improves the state-of-the-art sampling rate results for nonconvex matrix completion with no spurious local minima in \citet{ge2016matrix, ge2017no}. In addition, we investigated when the proposed nonconvex optimization results in accurate low-rank approximations even in presence of large condition numbers, large incoherence parameters, or rank mismatching. We also propose to apply the nonconvex optimization to memory-efficient Kernel PCA. Compared to the well-known Nystr\"{o}m methods, numerical experiments indicate that the proposed nonconvex optimization approach yields more stable results in both low-rank approximation and clustering.

\end{abstract}

\section{Introduction} 
\label{sec_intro}
Let $\mtx{M}$ be an $n \times n$ positive semidefinite matrix and let $r \ll n$ be a fixed integer. It is well known that a rank-$r$ approximation of $\mtx{M}$ can be obtained by truncating the spectral decomposition of $\mtx{M}$. To be specific, let $\mtx{M} = \sum_{i=1}^n \sigma_i \vct{u}_i \vct{u}_i^\top$ be the spectral decomposition with $\sigma_1 \geqslant \ldots \geqslant \sigma_n \geqslant 0$. Then, the best rank-$r$ approximation of $\mtx{M}$ is $\mtx{M}_r = \sum_{i=1}^r \sigma_i \vct{u}_i \vct{u}_i^\top$. If we denote $\boldsymbol{U}_r = [\sqrt{\sigma_1}\boldsymbol{u}_1\; \dots\; \sqrt{\sigma_r}\boldsymbol{u}_r]$, then the best rank-$r$ approximation of $\mtx{M}$ can be written as $\mtx{M} = \mtx{U}_r \mtx{U}_r^\top$. By the well-known Eckart-Young-Mirsky Theorem \citep{golub2012matrix}, $\mtx{U}_r$ is actually the global minimum (up to rotation) to the following nonconvex optimization:
\begin{equation*}
\label{eq:approximation_pos}
\begin{aligned}
\min_{\mtx{X} \in \mathbb{R}^{n \times r}} \| \mtx{X} \mtx{X}^\top - \mtx{M}\|_F^2.
\end{aligned}
\end{equation*} 
This factorization for low-rank approximation has been well-known in the literature \citep[see, e.g.,][]{burer2003nonlinear}.

This paper studies how to find a rank-$r$ approximation of $\mtx{M}$ in the case that only partial entries are observed. Let $\Omega\subset [n]\times [n]$ be a symmetric index set, and we assume that $\mtx{M}$ is only observed on the entries in $\Omega$. For convenience of discussion, this subsampling is represented as $\mathcal{P}_{\Omega}(\boldsymbol{M})$ in that $\mathcal{P}_{\Omega}(\boldsymbol{M})_{i,j} = M_{i,j}$ if $(i, j) \in \Omega$ and $\mathcal{P}_{\Omega}(\boldsymbol{M})_{i,j} = 0$ if $(i, j) \notin \Omega$. We are interested in the following question
\begin{center}
\emph{How to find a rank-$r$ approximation of $\mtx{M}$ in a scalable manner only through $\mathcal{P}_{\Omega}(\mtx{M})$?}
\end{center}

We propose to find such a low-rank approximation through the following nonconvex optimization, which has been exactly proposed in \citet{ge2016matrix, ge2017no} for matrix completion. Denote $\mtx{X} = \begin{bmatrix} \vct{x}_1, \ldots, \vct{x}_n\end{bmatrix}^\top \in \mathbb{R}^{n \times r}$. A rank-$r$ approximation of $\mtx{M}$ can be found through 
\begin{align}
\label{eq:obj_psd}
\min_{\boldsymbol{X}\in\mathbb{R}^{n\times r}} f(\boldsymbol{X}) & \coloneqq \frac{1}{2}\sum_{(i, j) \in \Omega} \left(\vct{x}_i^\top \vct{x}_j - M_{ij}\right)^2 +\lambda \sum_{i=1}^{n}[(\|\vct{x}_i\|_2-\alpha)_+]^4 \nonumber
\\
& \coloneqq \frac{1}{2}\|\mathcal{P}_{\Omega}(\boldsymbol{X}\boldsymbol{X}^\top-\boldsymbol{M})\|_F^2 +\lambda G_{\alpha}(\boldsymbol{X})
\end{align}
where $G_{\alpha}(\boldsymbol{X}) \coloneqq \sum_{i=1}^{n}[(\|\vct{x}_i\|_2-\alpha)_+]^4$. Following the framework of nonconvex optimization \emph{without initialization} in \citet{ge2016matrix, ge2017no}, our local-minimum based approximation for $\mtx{M}$ is 
$\mtx{M} \approx \widehat{\mtx{X}}\widehat{\mtx{X}}^\top$ where $\widehat{\mtx{X}}$ is \emph{any} local minimum of \prettyref{eq:obj_psd}.

Let's briefly discuss the memory and computational complexity to solve \eqref{eq:obj_psd} via gradient descent. If $\Omega$ is symmetric and does not contain the diagonal entries as later specified in Definition \ref{def:ber}, the updating rule of gradient decent 
\begin{equation}
\label{eq:gradient_descent}
\boldsymbol{X}^{(t+1)} = \boldsymbol{X}^{(t)}-\eta^{(t)}\nabla f(\boldsymbol{X}^{(t)})
\end{equation}
is equivalent to
\[
\vct{x}_i^{(t+1)}\coloneqq \vct{x}_i^{(t)} - \eta^{(t)}\left[2 \sum_{j: (i, j) \in \Omega} \left(\langle \vct{x}_i^{(t)}, \vct{x}_j^{(t)} \rangle - M_{i,j}\right)\vct{x}_j^{(t)} + \frac{4 \lambda}{\|\vct{x}_i^{(t)}\|_2} \left(\|\vct{x}_i^{(t)}\|_2 - \alpha \right)^3 1_{\{\|\vct{x}_i^{(t)}\|_2 \geqslant  \alpha\}} \vct{x}_i^{(t)}\right],
\]
where the memory cost is dominated by storing $\boldsymbol{X}^{(t)}$, $\boldsymbol{X}^{(t+1)}$, and $\mtx{M}$ on $\Omega$, which is generally $O(nr + |\Omega|)$. It is also obvious that the computational cost in each iteration is $O(|\Omega|r)$.

\subsection{Applications in memory-efficient kernel PCA}

Kernel PCA \citep{scholkopf1998nonlinear} is a widely used nonlinear dimension reduction technique in machine learning for the purpose of redundancy removal and preprocessing before prediction, classification or clustering. The method is implemented by finding a low-rank approximation of the kernel-based Gram matrix determined by the data sample. To be concrete, let $\vct{z}_1, \ldots, \vct{z}_n$ be a data sample of size $n$ and dimension $d$, and let $\mtx{M}$ be the $n \times n$ positive semidefinite kernel matrix determined by a predetermined kernel function $K(\vct{x}, \vct{y})$ in that $M_{ij} = K(\vct{z}_i, \vct{z}_j)$. Non-centered Kernel PCA with $r$ principal components amounts to finding the best rank-$r$ approximation of $\mtx{M}$. 

However, when the sample size is large, the storage of the kernel matrix itself becomes challenging. Consider the example when the dimension $d$ is in thousands while the sample size $n$ is in millions. The memory cost for the data matrix is $d \times n$ and thus in billions, while the memory cost for the kernel matrix $\mtx{M}$ is in trillions! On the other hand, if not storing $\mtx{M}$, the implementation of standard iterative algorithms of SVD will involve one pass of computing all entries of $\mtx{M}$ in each iteration, usually with formidable computational cost $O(n^2 d)$. A natural question arises: \emph{How to find low-rank approximations of $\mtx{M}$ memory-efficiently?}

The following two are among the most well-known memory-efficient Kernel PCA methods in the literature. One is Nystr\"{o}m method \citep{williams2001using}, which amounts to generating \emph{random partial columns} of the kernel matrix, then finding a low-rank approximation based on these columns. In order to generate random partial columns, uniform sampling without replacement is employed in \citet{williams2001using}, and different sampling strategies are proposed later \citep[e.g.,][]{drineas2005nystrom}. The method is convenient in implementation and efficient in both memory and computation, but relatively unstable in terms of approximation errors as will be shown in \prettyref{sec:simulations}. 

Another popular approach is stochastic approximation, e.g.,  Kernel Hebbian Algorithm (KHA) \citep{kim2005iterative}, which is memory-efficient and approaches the exact principal component solution as the number of iterations goes to infinity with appropriately chosen learning rate \citep{kim2005iterative}. However, based on our experience, the method usually requires careful tuning of learning rates even for very slow convergence.

It is also worth mentioning that the randomized one-pass algorithm discussed in, e.g., \citet{halko2011finding}, where the theoretical properties of a random-projection based low-rank approximation method were fully analyzed. However, although the one-pass algorithm does not require the storage of the whole matrix $\boldsymbol{M}$, in Kernel PCA one still needs to compute every entry of $\boldsymbol{M}$, which typically requires $O(n^2 d)$ computational complexity for kernel matrix.

As a result, we aim at finding a memory-efficient method as an alternative to the aforementioned approaches. In particular, we are interested  in a method with desirable empirical properties: memory-efficient, no requirement on one or multiple passes to compute the complete kernel matrix, no requirement to tune the parameters carefully, and yielding stable results. To this end, we propose the following method based on entries sampling and nonconvex optimization: In the first step, $\Omega$ is generated to follow an Erd\H{o}s-R\'{e}nyi random graph with parameter $p$ later specified in Definition \ref{def:ber}, and then a partial kernel matrix $\mathcal{P}_{\Omega}(\mtx{M})$ is generated in that $M_{i,j} = K(\vct{z}_i, \vct{z}_j)$ for $(i, j) \in \Omega$. In the second step, the nonconvex optimization \eqref{eq:obj_psd} is implemented through gradient descent \eqref{eq:gradient_descent}. Any local minimum of \prettyref{eq:obj_psd}, $\widehat{X}$, is a solution of approximate kernel PCA in that $\mtx{M} \approx \widehat{X}\widehat{X}^\top$.


To store the index set $\Omega$ and the sampled entries of $\mtx{M}$ on $\Omega$, the memory cost in the first step is $O(|\Omega|)$, which is comparable to the memory cost $O(nr + |\Omega|)$ in the second step. As to the computational complexity, besides the generation of $\Omega$, the computational cost in the first step is typically $O(|\Omega|d)$, e.g., when the radial kernels or polynomial kernels are employed. This could be dominating the per-iteration computational complexity $O(|\Omega|r)$ in the second step when the target rank $r$ is much smaller than the original dimension $d$.

Partial entries sampling plus nonconvex optimization has been proposed in the literature for scalable robust PCA and matrix completion \citep{yi2016fast}. However, to the best of our knowledge, our work is the first to apply such an idea to memory-efficient Kernel PCA. Moreover, the underlying signal matrix is assumed to be exactly low-rank in \citet{yi2016fast} while we make no assumptions on the positive semidefinite kernel matrix $\mtx{M}$. Entry-sampling has been proposed in \citet{achlioptas2002sampling, achlioptas2007fast} for scalable low-rank approximation. In particular, it is used to speed up Kernel PCA in \citet{achlioptas2002sampling}, but spectral methods are subsequently employed after entries sampling as opposed to nonconvex optimization. Empirical comparisons between spectral methods and nonconvex optimization will be demonstrated in \prettyref{sec:simulations}. It is also noteworthy that matrix completion techniques have been applied to certain kernel matrices when it is costly to generate each single entry \citep{Graepel2002kernel, Paisley2010kernel}, wherein the proposed methods are not memory-efficient. In contrast, our method is memory-efficient in order to serve a different purpose.

\subsection{Related work and our contributions}
In recent years, a series of papers have been proposed to study nonconvex matrix completion \citep[see, e.g.,][]{rennie2005fast, keshavan2010matrix, keshavan2010JMLR,jain2013low, zhao2015nonconvex, sun2015guaranteed, chen2015fast, yi2016fast, zheng2016convergence, ge2016matrix, ge2017no}. Interested readers are referred to \citet{balcan2017optimal}, where required sampling rates in these papers are summarized in Table 1 therein. Compared to convex approaches for matrix completion \citep[e.g.,][]{candes2009exact}, these nonconvex approaches are not only more computationally efficient, but also more convenient in storing. For the same reason, nonconvex optimization approaches have also been investigated for other low-rank recovery problems including phase retirval \citep[e.g.,][]{candes2015phase,sun2016geometric,cai2016optimal}, matrix sensing \citep[e.g.,][]{zheng2015convergent,tu2015low}, blind deconvolution \citep[e.g.,][]{li2018rapid}, etc. 

Our present work follows the framework of local minimum analysis for nonconvex optimization in the literature. For example, \citet{baldi1989neural} has described the nonconvex landscape of the quadratic loss for PCA. \citet{JMLR:v16:loh15a} studies the local minima of regularized M-estimators. \citet{sun2016geometric} studies the global geometry of the phase retrieval problem. The conditions for no spurious local minima have been investigated in \citet{BNS16} and \citet{ge2016matrix} for nonconvex matrix sensing and completion, respectively. The global geometry of nonconvex objective functions with underlying symmetric structures, including low-rank symmetric matrix factorization and sensing, has been studied in \citet{li2016symmetry}. Global geometry of rectangular matrix factorization and sensing is studied \citet{zhu2017global}, where the issues of under-parameterization and over-parameterization have been investigated. Similar analysis is extend to general low-rank optimization problems in \citet{li2017geometry}. Matrix factorization is further studied in \citet{jin2017escape} with a novel geometric characterization of saddle points, and this idea is later extended in \citet{ge2017no}, where a unified geometric analysis framework is proposed to study the landscapes of nonconvex matrix sensing, matrix completion and robust PCA.



Among these results, \citet{ge2016matrix} and \citet{ge2017no} are highly relevant to our work in both methodological and technical terms. In fact, exactly the same nonconvex optimization problem \eqref{eq:obj_psd} has been studied in \citet{ge2016matrix, ge2017no} for matrix completion from missing data. To be specific, these papers show that any local minimum $\widehat{\mtx{X}}$ yields $\mtx{M} = \widehat{\mtx{X}}\widehat{\mtx{X}}^\top$, as long as $\mtx{M}$ is exactly rank-$r$, the condition number $\kappa_r \coloneqq  \sigma_1/\sigma_r$ is well-bounded, the incoherence parameter of the eigenspace of $\mtx{M}$ is well-bounded, and the sampling rate is greater than a function of these quantities. The case with additive stochastic noise has also been discussed in \citet{ge2016matrix}.

In contrast, our paper studies the theoretical properties of $\widehat{\mtx{X}}\widehat{\mtx{X}}^\top$ with \emph{no assumptions} on $\mtx{M}$. There are actually two questions of interest: how close $\widehat{\mtx{X}}\widehat{\mtx{X}}^\top$ is from $\mtx{M}$, and how close $\widehat{\mtx{X}}\widehat{\mtx{X}}^\top$ is from $\mtx{M}_r$ (recall that $\mtx{M}_r$ is the best rank-$r$ approximation of $\mtx{M}$ by spectral truncation). In comparison to \citet{ge2016matrix,ge2017no}, our main contributions to be introduced in the next section include the following:
\begin{itemize}
\item Our main result \prettyref{thm:main_psd} that characterizes how well any local-minimum based rank-$r$ factorization $\widehat{\mtx{X}}\widehat{\mtx{X}}^\top$ approximates $\mtx{M}$ or $\mtx{M}_r$ requires no assumptions imposed on $\mtx{M}$ regarding its rank, eigenvalues and eigenvectors. The sampling rate is only required to satisfy $p\geqslant  C(\log n /n)$ for some absolute constant $C$. Therefore, for applications such as memory-efficient Kernel PCA, our framework provides more suitable guidelines than \citet{ge2016matrix, ge2017no}. In fact, kernel matrices are in general of full rank and their condition numbers and incoherence parameters may not satisfy the strong assumptions in \citet{ge2016matrix, ge2017no}.

\item When $\mtx{M}$ is assumed to be exactly low-rank as in \citet{ge2016matrix, ge2017no}, Corollary \ref{coro_ex_psd_mu} improves the state-of-the-art no-spurious-local-minima results in \citet{ge2016matrix, ge2017no} for exact nonconvex matrix completion in terms of sampling rates. To be specific, assuming both condition numbers and incoherence parameters are on the order of $O(1)$, our result improves the result in \cite{ge2017no} from $\widetilde{O}(r^4/n)$ to $\widetilde{O}(r^2/n)$.

\item \prettyref{thm:main_psd} also implies the conditions under which the nonconvex optimization \eqref{eq:obj_psd} yields good low-rank approximation of $\mtx{M}$ in the cases of large condition numbers, high incoherence parameters, or rank-mismatching.
\end{itemize}

\begin{figure}[ht]
\center
\begin{tikzpicture}[line join=round]

	 \draw (-7,0) -- (7,0)  ;
 \draw [blue,line width=0.45mm] plot [smooth, tension=0.8] coordinates { (-4,-2) (-2,3) (0,0)  (2,2)(4,-2)}; 
	\draw [green,line width=0.45mm] plot [smooth, tension=0.6] coordinates { (-4,-5)  (-2,-1)(-1.5,-1.5)(-1,-0.25)(-0.5,-1) (-0.2,-0.75)(0.5,-1)(1,-0.5)(1.5,-1.5)(2,-1) (4,-5)};
	\filldraw (2,2) node[above]{$K ( \boldsymbol{X})$}; 
	 \filldraw (1,-0.5) node[right]{$-f(\boldsymbol{X})$};
	
	 \filldraw (0,0) circle(2pt) node[below ]{$\boldsymbol{U}_r$}; 
	 
	 \draw [line width=0.05mm] (-1.85,-1) -- (-1.85,-3);
	 \draw [line width=0.05mm] (1.85,-1) -- (1.85,-3);
	 \draw [<-] (-1.85,-2.5)--(-1.25,-2.5);
		\draw [->] (1.25,-2.5)--(1.85,-2.5);
		 \filldraw (0,-2.5) node[align=center]{ span of\\ local minima \\ of $f(\boldsymbol{X})$};

	 \draw [line width=0.05mm] (-3.4 ,-0.25) -- (-3.4, -4);
	 \draw  [line width=0.05mm](3.25,-0.25) -- (3.25,-4);
	 \draw [<-] (-3.4,-3.75)--(-2.75,-3.75);
		\draw [->] (2.75,-3.75)--(3.25,-3.75);
		 \filldraw (0,-3.75) node[align=center]{ span of $\{\boldsymbol{X}\in\mathbb{R}^{n\times r} \mid K(\boldsymbol{X})\geqslant  0\}$};

\end{tikzpicture}
\captionsetup{justification=centering,margin=2cm}
\caption{ Landscape of $-f(\boldsymbol{X}), K(\boldsymbol{X})$ and $\boldsymbol{U}_r$.}
\label{fg01}
\end{figure}
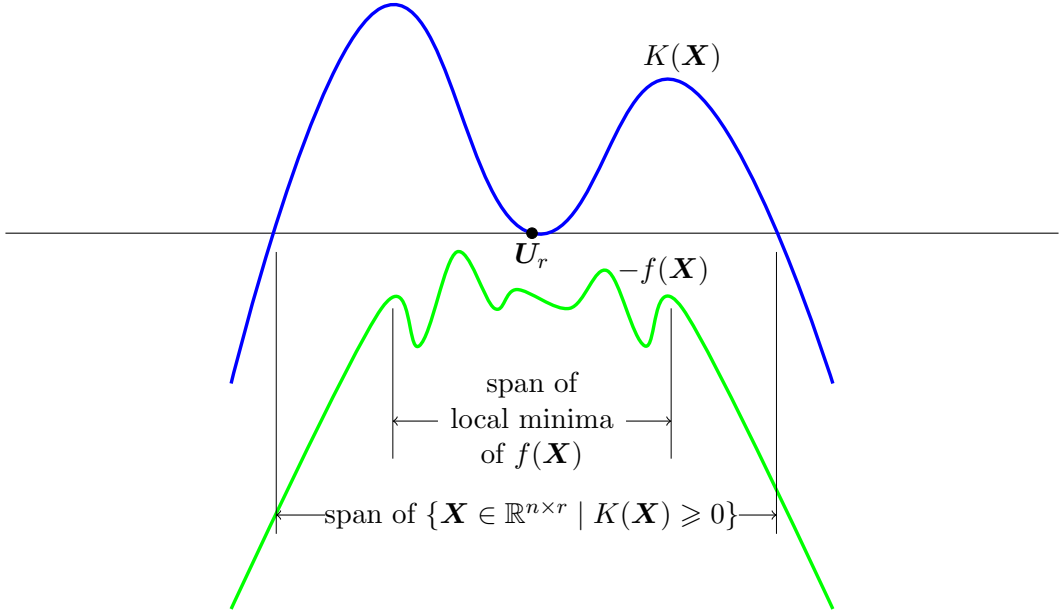

On the other hand, our paper benefits from \citet{ge2016matrix, ge2017no} in various aspects. In order to characterize the properties of any local minimum $\widehat{\mtx{X}}$, we follow the idea in \citet{ge2017no} to combine the first and second order conditions of local minima linearly to construct an auxiliary function, denoted as $K(\mtx{X})$ in our paper, and consequently all local minima satisfy the inequality $K(\widehat{\mtx{X}}) \geqslant  0$ as illustrated in Figure \ref{fg01}. If $\mtx{M}$ is exactly rank-$r$ and its eigenvalues and eigenvectors satisfy particular properties, \citet{ge2017no} shows that $K(\mtx{X}) \leqslant 0$  for all $\mtx{X}$ as long as the sampling rate is large enough. This argument can be employed to prove that there is no spurious local minima.

However, $K(\mtx{X}) \leqslant 0$ is not always true if no assumptions are imposed on $\mtx{M}$, so we instead focus on analyzing the inequality $K(\widehat{\mtx{X}}) \geqslant  0$ directly in a model-free setup. Among a few novel technical ideas, it is worth highlighting the deterministic inequality (Lemma \ref{lemma_concern}) that controls the difference between the function $K(\mtx{X})$ and its population version $\E[K(\mtx{X})]$ in a uniform and model-free manner.

\subsection{Organization and notations}
The remainder of the paper is organized as follows: Our main theoretical results are stated in \prettyref{sec:theory}; Numerical simulations and applications in memory-efficient KPCA are given in \prettyref{sec:simulations}. Proofs are deferred to \prettyref{sec:proofs}. 

We use bold letters to denote matrices and vectors. For any vectors $\vct{u}$ and $\boldsymbol{v}$, $\|\boldsymbol{u}\|_2$ denotes its $\ell_2$ norm, and $\langle \boldsymbol{u},\boldsymbol{v}\rangle$ their inner product. For any matrix $\boldsymbol{M}\in\mathbb{R}^{n \times n}$, $M_{i,j}$ denotes its $(i,j)$-th entry, $\boldsymbol{M}_{i,\cdot} = (M_{i,1},M_{i,2},\dots,M_{i,n})^\top$ its $i$-th row of $\boldsymbol{M}$, and $\boldsymbol{M}_{\cdot,j} = (M_{1,j},M_{2,j},\dots,M_{n,j})^\top$ its $j$-th column. Moreover, we use $\|\boldsymbol{M}\|$, $\|\boldsymbol{M}\|_*$, $\|\boldsymbol{M}\|_F$, $\|\boldsymbol{M}\|_{\ell_\infty}\coloneqq \max_{i,j}|M_{i,j}|$, $\|\boldsymbol{M}\|_{2,\infty}\coloneqq\max_{i}\|\boldsymbol{M}_{i,\cdot}\|_2$ to denote its spectral norm, nuclear norm, Frobenius norm, elementwise max norm and $\ell_{2,\infty}$ norm, respectively. The vectorization of $\boldsymbol{M}$ is represented by $\operatorname{vec}(\boldsymbol{M}) = (M_{1,1},M_{2,1},\dots,M_{1,2},\dots, M_{n,n})^\top$. For matrices $\boldsymbol{M},\boldsymbol{N}$ of the same size, denote $\langle \boldsymbol{M},\boldsymbol{N} \rangle = \sum_{i,j} M_{i,j}N_{i,j}= \text{trace}\left(\boldsymbol{M}^\top \boldsymbol{N}\right)$. Denote by $\nabla f(\boldsymbol{M})\in\mathbb{R}^{n\times n}$ and $\nabla^2 f(\boldsymbol{M})\in\mathbb{R}^{n^2\times n^2}$ the gradient and Hessian of $f(\mtx{M})$. 

Denote $[x]_+ = \max\{x,0\}$. We use $\boldsymbol{J}$ to denote a matrix whose all entries equal to one. We use $C,C_1,C_2,\dots$ to denote absolute constants, whose values may change from line to line.

\section{Model-free approximation theory}
\label{sec:theory}

\subsection{Main results}
The following sampling scheme is employed throughout the paper:
\begin{definition}[Off-diagonal symmetric independent $\text{Ber}(p)$ model] 
\label{def:ber}
Assume the index set $\Omega$ consists only of off-diagonal entries that are sampled symmetrically and independently with probability $p$, i.e.,
\begin{enumerate}
\item $(i, i) \notin \Omega$ for all $i =1, \dots, n$;
\item For all $i<j$, sample $(i,j)\in\Omega$ independently with probability $p$; 
\item For all $i>j$, $(i,j)\in \Omega$ if and only if $(j,i)\in \Omega$.
\end{enumerate}
\end{definition}
Here we assume all diagonal entries are not in $\Omega$ for the generality of the formulation, although they are likely to be obtained in practice. For instance, all diagonal entries of the radial kernel matrix are ones. For any index set $\Omega \subset [n] \times [n]$, define the associated $0$-$1$ matrix $\boldsymbol{\Omega}\in \{0,1\}^{n\times n}$ such that $\Omega_{i,j}=1$ if and only if $(i, j) \in \Omega$. Then we can write $\mathcal{P}_{\Omega}(\boldsymbol{X}) = \boldsymbol{X}\circ \boldsymbol{\Omega}$ where $\circ$ is the Hadamard product.  

Assume that the positive semidefinite matrix $\boldsymbol{M}$ has the spectral decomposition
\begin{equation}
\label{eq:spectral}
 \boldsymbol{M} = \sum_{i=1}^r \sigma_i \boldsymbol{u}_i\boldsymbol{u}_i^\top + \sum_{i=r+1}^n \sigma_i \boldsymbol{u}_i\boldsymbol{u}_i^\top  := \boldsymbol{M}_r+\boldsymbol{N},
\end{equation}
where $\sigma_1\geqslant  \sigma_2\geqslant  \dots\geqslant  \sigma_n \geqslant  0$ are the spectrum, $\boldsymbol{u}_i\in \mathbb{R}^{n }$ are unit and mutually perpendicular eigenvectors. The matrix $\mtx{M}_r \coloneqq \sum_{i=1}^r \sigma_i \boldsymbol{u}_i\boldsymbol{u}_i^\top$ is the best rank-$r$ approximation of $\boldsymbol{M}$ and $\boldsymbol{N}\coloneqq \sum_{i=r+1}^n \sigma_i \boldsymbol{u}_i\boldsymbol{u}_i^\top$ denotes the residual part. In the case of multiple eigenvalues, the order in the eigenvalue decomposition \eqref{eq:spectral} may not be unique. In this case, we consider the problem for any fixed order in \eqref{eq:spectral} with the fixed $\mtx{M}_r$.

\begin{theorem} 
\label{thm:main_psd}
Let $\boldsymbol{M} \in\mathbb{R}^{n\times n}$ be a positive semidefinite matrix with the spectral decomposition \prettyref{eq:spectral}. Let $\Omega$ be sampled according to the off-diagonal symmetric $\text{Ber}(p)$ model with $ p\geqslant C_1\frac{\log n}{n}$ for some absolute constant $C_1$. Then in an event $E$ with probability $\mathbb{P}[E]\geqslant  1-2n^{-3}$, as long as the tuning parameters $\alpha$ and $\lambda$ satisfy $100\sqrt{\|\boldsymbol{M}_r\|_{\ell_\infty}}\leqslant\alpha\leqslant 200\sqrt{\|\boldsymbol{M}_r\|_{\ell_\infty}}$ and $100 \|\boldsymbol{\Omega}-p\boldsymbol{J}\| \leqslant \lambda \leqslant 200 \|\boldsymbol{\Omega}-p\boldsymbol{J}\| $, any local minimum $\widehat{\boldsymbol{X}}\in\mathbb{R}^{n\times r}$ of \eqref{eq:obj_psd} satisfies
\begin{equation} \label{eq_040} 
\begin{split}
	 \left\|\widehat{\boldsymbol{X}}\widehat{\boldsymbol{X}}^\top -\boldsymbol{M}_r \right\|_F^2  \leqslant &  C_2  \sum_{i=1}^r \left\{\left[C_3  \left(\sqrt{\frac{n }{p}} +\frac{\log n}{p} \right)\|\boldsymbol{M}_r\|_{\ell_\infty}  + C_3 \sigma_{2r+1-i} -\sigma_i \right]_+\right\}^2 \\
	 & +  C_2\frac{[p(1-p)n+\log n]r \|\boldsymbol{N}\|_{\ell_\infty}^2}{p^2} 
\end{split}
\end{equation}
and
\begin{equation} \label{eq_062}
\begin{split}
	\left\|\widehat{\boldsymbol{X}}\widehat{\boldsymbol{X}}^\top - \boldsymbol{M}\right\|_F^2  \leqslant & C_2    \sum_{i=1}^r \left\{\left[C_3  \left(\sqrt{\frac{n }{p}} +\frac{\log n}{p} \right)\|\boldsymbol{M}_r\|_{\ell_\infty}   + C_3 \sigma_{2r+1-i} -\sigma_i \right]_+\right\}^2  \\
	&+  C_2 \frac{[p(1-p)n+\log n]r \|\boldsymbol{N}\|_{\ell_\infty}^2}{p^2}   + \|\boldsymbol{N}\|_F^2 
\end{split}
\end{equation}
with $C_2,C_3$ absolute constants.


\end{theorem}

Model-free low-rank approximation from partial entries has been studied for for spectral estimators in the literature. For example, under the settings of Theorem \ref{thm:main_psd}, the spectral low-rank approximation (denoted as $\boldsymbol{M}_{\textrm{approx}}$) discussed in \citet[Theorem 1.1]{keshavan2010JMLR} is guaranteed to satisfy
\begin{equation*} 
 \|\boldsymbol{M}_{\textrm{approx}}- \boldsymbol{M}_r\|_F^2 \leqslant C\left\{ \frac{nr\|\boldsymbol{M}_r\|_{\ell_\infty}^2}{p} + \frac{r\|\mathcal{P}_{\Omega}(\boldsymbol{N})\|^2}{p^2} \right\},
\end{equation*}
with high probability. However, this cannot imply any exact recovery results even when $\mtx{M}$ is of low rank and the sampling rate $p$ satisfies the conditions specified in \citet{ge2017no}. Similarly, the SVD-based USVT estimator introduced in \citet{chatterjee2015matrix} does not imply exact recovery. In contrast, as will be discussed in the next subsection, \prettyref{thm:main_psd} implies that any local minimum of \eqref{eq:obj_psd} yields exact recovery of $\mtx{M}$ with high probability under milder conditions than those in \citet{ge2017no}.

\subsection{Implications in exact matrix completion}
\label{sec:exact_mc}
Assume in this subsection that the positive semidefinite matrix $\boldsymbol{M}$ is exactly rank-$r$, i.e., 
\begin{equation}
\label{eq:decomp_low_rank}
\mtx{M} = \mtx{M}_r = \sum_{i=1}^r \sigma_i \boldsymbol{u}_i\boldsymbol{u}_i^\top = \boldsymbol{U}_r\boldsymbol{U}_r^\top
\end{equation}
where $\boldsymbol{U}_r = [\sqrt{\sigma_1}\boldsymbol{u}_1\; \dots\; \sqrt{\sigma_r}\boldsymbol{u}_r]$. Furthermore, we assume its condition number 
$\kappa_r = \frac{\sigma_1}{\sigma_r}$ and eigen-space incoherence parameter $\mu_r = \frac{n}{r}\max_i \sum_{j=1}^r u_{i,j}^2$ \citep{candes2009exact} are well-bounded. This is a standard setup in the literature of nonconvex matrix completion \citep[e.g.,][]{keshavan2010matrix, sun2015guaranteed, chen2015fast, zheng2016convergence, ge2016matrix, yi2016fast, ge2017no}.

Notice that \citet{ge2016matrix} introduces a slightly different version of incoherence
\begin{equation}
\label{eq:rong_spikiness}
\widetilde{\mu}_r \coloneqq \frac{\sqrt{n} \|\boldsymbol{U}_r\|_{2,\infty}}{\|\boldsymbol{U}_r\|_F}= \sqrt{\frac{n\|\boldsymbol{M}_r\|_{\ell_\infty}}{\operatorname{trace}(\boldsymbol{M}_r)}}
\end{equation}
as a measure of spikiness. (Note that this is different from the spikiness defined in \citet{negahban2012restricted}.) By $\|\boldsymbol{M}_r\|_{\ell_\infty} = \| \boldsymbol{U}_r \|_{2,\infty}^2 =  \max_i \sum_{j=1}^r \sigma_j u_{i,j}^2$, the following relationship between $\mu$ and $\widetilde{\mu}$ is straightforward
\begin{equation}
\label{eq_100}
 \frac{\widetilde{\mu}_r^2}{\kappa_r} \leqslant \frac{\widetilde{\mu}_r^2 \operatorname{trace}(\boldsymbol{M}_r)}{r\sigma_1} =\frac{n\|\boldsymbol{M}_r\|_{\ell_\infty}}{r\sigma_1} \leqslant \mu_r  \leqslant \frac{n\|\boldsymbol{M}_r\|_{\ell_\infty}}{r\sigma_r} = \frac{\widetilde{\mu}_r^2 \operatorname{trace}(\boldsymbol{M}_r)}{r\sigma_r}\leqslant \kappa_r \widetilde{\mu}_r^2.
\end{equation}

By the fact $\|\mtx{M}\|_{\ell_\infty} \leqslant \frac{r}{n} \sigma_1 \mu_r$, \prettyref{thm:main_psd} implies the following exact low-rank recovery results:

\begin{corollary}
\label{coro_ex_psd_mu}
Under the assumptions of Theorem \ref{thm:main_psd}, if we further assume $\rank(\boldsymbol{M})=r$ (i.e., $\mtx{M} = \mtx{M}_r$) and
\begin{equation*}  
 p \geqslant C \max\left\{\frac{\mu_r  r \kappa_r \log n}{n}, \frac{\mu_r^2  r^2\kappa_r^2 }{n}\right\}   
\end{equation*}
or
\begin{equation*}
 p \geqslant  C \max\left\{ \frac{\widetilde{\mu}_r^2 r \kappa_r \log n}{n}, \frac{\widetilde{\mu}_r^4 r^2\kappa_r^2}{n} \right\}  
\end{equation*}
for some absolute constant $C$, then in an event $E$ with probability $\mathbb{P}[E]\geqslant  1-2n^{-3}$, any local minimum $\widehat{\boldsymbol{X}}\in\mathbb{R}^{n\times r}$ of objective function $f(\boldsymbol{X})$ defined in \eqref{eq:obj_psd} satisfies $\widehat{\boldsymbol{X}}\widehat{\boldsymbol{X}}^\top  = \boldsymbol{M}$.  
\end{corollary}
The proof is straightforward and deferred to the appendix. Notice that our results are better than the state-of-the-art results for no spurious local minimum in \citet{ge2017no}, where the required sampling rate is $p\geqslant  \frac{C}{n}\mu_r^3 r^4\kappa_r^4\log n$ (which also implies $p\geqslant  \frac{C}{n}\widetilde{\mu}_r^6 r^4\kappa_r^{7}\log n$ by \eqref{eq_100}).

\subsection{Examples}
Besides improving the state-of-the-art no-spurious-local-minima results in nonconvex matrix completion, Theorem \ref{thm:main_psd} is also capable of explaining some nontrivial phenomena in low-rank matrix completion in the presence of large condition numbers, high incoherence parameter, or mismatching between the selected and true ranks.

\subsubsection{Nonconvex matrix completion with large condition numbers and high eigen-space incoherence parameters}\label{sec:exmp1}
Assume here $\mtx{M}$ is exactly rank-$r$ and its spectral decomposition is denoted as in \eqref{eq:decomp_low_rank}. However, we assume that $\mu_r$ and $\kappa_r$ can be extremely large, while the condition number and incoherence parameter for $\mtx{M}_{r-1} = \sum_{i=1}^{r-1} \sigma_i \boldsymbol{u}_i\boldsymbol{u}_i^\top$, i.e., $\kappa_{r-1} = \frac{\sigma_1}{\sigma_{r-1}}$ and $\mu_{r-1} = \frac{n}{r-1}\max_{i}\sum_{j=1}^{r-1} u_{i,j}^2$, are well-bounded. We are interested in figuring out when the local minimum based rank-$r$ factorization $\widehat{\mtx{X}}\widehat{\mtx{X}}^\top$ approximates the original $\mtx{M}$ well.


By $\|\boldsymbol{M}_r\|_{\ell_\infty} =  \max_i \sum_{j=1}^r \sigma_j u_{i,j}^2$, we have
\[
\|\boldsymbol{M}_r\|_{\ell_\infty}\leqslant \frac{r-1}{n} \sigma_1 \mu_{r-1}+\sigma_r\|\boldsymbol{u}_r\|_{\infty}^2.
\]
Then by Theorem \ref{thm:main_psd}, if 
\[
 p\geqslant  C\max\left\{ \frac{\left[\mu_{r-1} \kappa_{r-1} (r-1)+n\frac{\sigma_r}{\sigma_{r-1}} \|\boldsymbol{u}_r\|_{\infty}^2\right]\log n}{ n}, \frac{\left[\mu_{r-1} \kappa_{r-1} (r-1)+n\frac{\sigma_r}{\sigma_{r-1}} \|\boldsymbol{u}_r\|_{\infty}^2\right]^2}{ n}\right\}	
\]
with some absolute constant $C$, in an event $E$ with probability $\mathbb{P}[E]\geqslant  1-2n^{-3}$, for any local minimum $\widehat{\boldsymbol{X}}\in\mathbb{R}^{n\times r}$ of \eqref{eq:obj_psd}, $\|\widehat{\boldsymbol{X}}\widehat{\boldsymbol{X}}^\top - \boldsymbol{M}\|_F^2 \leqslant \frac{1}{100}\sigma_{r-1}^2$ holds. In other words, the relative approximation error satisfies $RE \coloneqq \frac{\|\widehat{\boldsymbol{X}}\widehat{\boldsymbol{X}}^\top - \boldsymbol{M}\|_F}{\|\boldsymbol{M}\|_F}\leqslant \frac{1}{10\sqrt{r-1}}$.

Notice that $\|\boldsymbol{u}_r\|_{\infty}^2 \leqslant \frac{r}{n} \mu_r$ and $\frac{\sigma_r}{\sigma_{r-1}} = \frac{\kappa_{r-1}}{\kappa_{r}}$, so the above sampling rate requirement is satisfied as long as  $\frac{\mu_r}{\kappa_r} \leqslant C \mu_{r-1}$ and
\[
		p\geqslant C\max\left\{\frac{\mu_{r-1} \kappa_{r-1} r\log n}{n},\frac{\mu_{r-1}^2 \kappa_{r-1}^2 r^2}{n}\right\}.
\]

%
%


\subsubsection{Rank mismatching}
In this subsection, $\boldsymbol{M}$ is assumed to be exactly rank-$R$, i.e., 
\begin{equation*}
\mtx{M} = \mtx{M}_R = \sum_{i=1}^R \sigma_i \boldsymbol{u}_i\boldsymbol{u}_i^\top = \boldsymbol{U}_R\boldsymbol{U}_R^\top
\end{equation*}
where $\boldsymbol{U}_R = [\sqrt{\sigma_1}\boldsymbol{u}_1\; \dots\; \sqrt{\sigma_R}\boldsymbol{u}_R]$. However, we consider the case that the selected rank $r$ is not the same as the true rank $R$, i.e., rank mismatching. As with \prettyref{sec:exact_mc}, we assume the condition number $\kappa_R = \frac{\sigma_1}{\sigma_R}$ and eigen-space incoherence parameter $\mu_R = \frac{n}{R}\max_i \sum_{j=1}^R \sigma_j u_{i,j}^2$ are well-bounded. As with \eqref{eq_100}, there holds $\|\mtx{M}\|_{\ell_\infty} \leqslant \frac{R}{n} \sigma_1 \mu_R$.

\paragraph{Case 1: $R<r$.} 
Theorem \ref{thm:main_psd} implies that if 
\[
 p\geqslant C\max\left\{ \frac{\mu_R \kappa_R R \log n}{ n},\frac{\mu_R^2 \kappa_R^2 R^2}{n}\right\}  
\]
for some absolute constant $C$, then in an event $E$ with probability $\mathbb{P}[E]\geqslant  1-2n^{-3}$, any local minimum $\widehat{\boldsymbol{X}}\in \mathbb{R}^{n\times r}$ of \eqref{eq:obj_psd} yields $\|\widehat{\boldsymbol{X}}\widehat{\boldsymbol{X}}^\top - \boldsymbol{M}\|_F^2 \leqslant \frac{1}{100}(r-R)\sigma_R^2$. This further yields the relative approximation error bound $RE \coloneqq \frac{\|\widehat{\boldsymbol{X}}\widehat{\boldsymbol{X}}^\top - \boldsymbol{M}\|_F}{\|\boldsymbol{M}\|_F}\leqslant \frac{1}{10}\sqrt{\frac{r-R}{R}}$.

\paragraph{Case 2: $R>r$.}
Recall that $\|\mtx{M}_r\|_{\ell_\infty} \leqslant \frac{r}{n} \sigma_1 \mu_r$. Moreover, 
\[
\|\boldsymbol{N}\|_{\ell_\infty} = \max_i \sum_{j=r+1}^R \sigma_j u_{i,j}^2 \leqslant \sigma_{r+1} \left(\max_i \sum_{j=1}^R  u_{i,j}^2\right) = \frac{\mu_R R}{n}\sigma_{r+1}.
\]
Theorem \ref{thm:main_psd} implies that if 
\[
 p\geqslant C\max\left\{\frac{\mu_r  r \kappa_r \log n}{n}, \frac{\mu_r^2  r^2\kappa_r^2 }{n}, \frac{\mu_R^2 R^3}{n}\right\}  
\] 
for some absolute constant $C$, then with high probability, any local minimum $\widehat{\boldsymbol{X}}\in \mathbb{R}^{n\times r}$ of \eqref{eq:obj_psd} yields
\[
\|\widehat{\boldsymbol{X}}\widehat{\boldsymbol{X}}^\top - \boldsymbol{M}_r\|_F^2   \leqslant C(\sigma_{r+1}^2 + \ldots + \sigma_{2r}^2),
\]
which implies that the relative error is well-controlled as long as $\sigma_{r+1}^2 + \ldots + \sigma_R^2$ accounts for a small proportion in $\sigma_1^2 + \ldots + \sigma_R^2$.

If we assume that $ 2C_3 \sigma_{r+1}< \sigma_r$ where $C_3$ is specified in Theorem \ref{thm:main_psd}, under the same sampling rate requirement as above, Theorem \ref{thm:main_psd} implies a much sharper result:
\[
\|\widehat{\boldsymbol{X}}\widehat{\boldsymbol{X}}^\top - \boldsymbol{M}_r\|_F^2 \leqslant  \frac{1}{100} \sigma_{r+1}^2,
\]
which yields the following (perhaps surprising) relative approximation error bound 
\[
RE \coloneqq \frac{\|\widehat{\boldsymbol{X}}\widehat{\boldsymbol{X}}^\top - \boldsymbol{M}_r\|_F}{\|\boldsymbol{M}_r\|_F}\leqslant \frac{1}{10}\sqrt{\frac{\sigma_{r+1}^2}{\sigma_1^2 + \ldots + \sigma_r^2}} \leqslant \frac{1}{10\sqrt{r}}.
\]

\section{Experiments}
\label{sec:simulations}
In the following simulations where the nonconvex optimization \eqref{eq:obj_psd} is solved, the initialization $\boldsymbol{X}^{(0)}$ is constructed randomly with i.i.d. normal entries with mean $0$ and variance $1$. The step size $\eta^{(t)}$ for the gradient descent \eqref{eq:gradient_descent} is determined by Armijo's rule \citep{armijo1966minimization}. The gradient descent algorithm is implemented with sparse matrix storage in Section \ref{sec:KPCA_num} for the purpose of memory-efficient KPCA, while with full matrix storage in Section \ref{sec:lowrank_num} to test the performance of general low-rank approximations from missing data. In each experiment, the iterations will be terminated when $\|\nabla f(\boldsymbol{X}^{(t)})\|_F\leqslant 10^{-3}$ or $\|\eta^{(t)}\nabla f(\boldsymbol{X}^{(t)})\|_F\leqslant 10^{-10}$ or the number of iterations surpasses $10^3$. All methods are implemented in MATLAB. The experiments are running on a virtual computer with Linux KVM, with 12 cores of 2.00GHz Intel Xeon E5 processor and 16 GB memory.

\subsection{Numerical simulations}
\label{sec:lowrank_num}
In this section, we conduct numerical tests on the nonconvex optimization \eqref{eq:obj_psd} under different settings of spectrum for the $500\times 500$ positive semidefinite matrix $\mtx{M}$, whose eigenvectors are the same as the left singular vectors of a random $500\times 500$ matrix with i.i.d. standard normal entries. The generation of eigenvalues for $\mtx{M}$ will be further specified in each test. For each generated $\boldsymbol{M}$, the nonconvex optimization \eqref{eq:obj_psd} is implemented for 50 times with independent $\Omega$'s generated under the off-diagonal symmetric independent Ber($p$) model. To implement the gradient descent algorithm \eqref{eq:gradient_descent}, set $\alpha = 100\|\boldsymbol{M}\|_{\ell_\infty}$ and $\lambda = 100\|\boldsymbol{\Omega}-p\boldsymbol{J}\|$ (the performances of our method are empirically not sensitive to the choices of the tuning parameters). In each single numerical experiment, we also conduct spectral method proposed in \citet{achlioptas2002sampling} to obtain an approximate low-rank approximation of $\mtx{M}$ for the purpose of comparison.

\subsubsection{Full rank case}\label{small02_sim}
\begin{figure}[!ht]
 \begin{subfigure}[t]{0.5\textwidth}
	 \includegraphics[width=\textwidth]{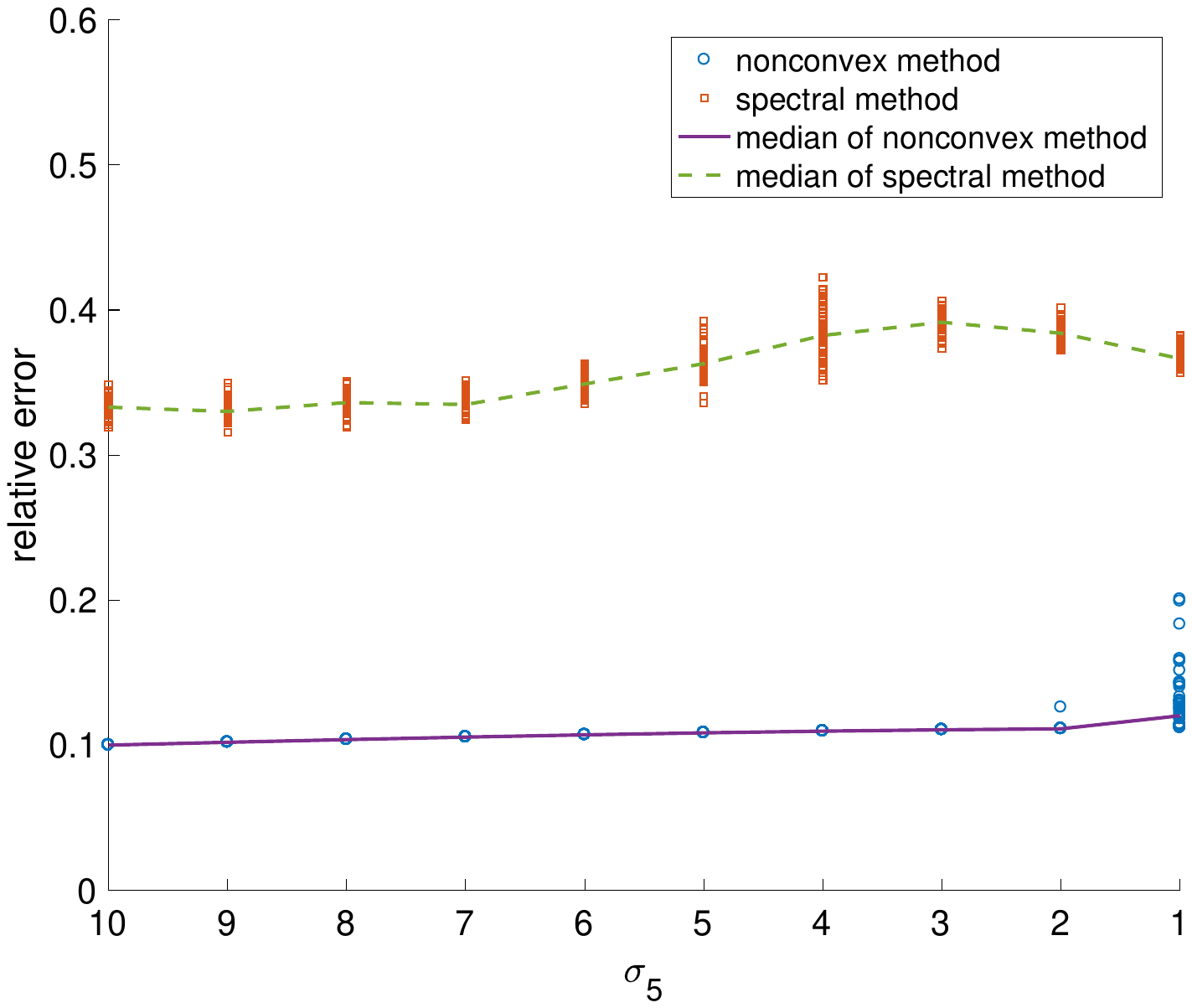}
	 \caption{   Relative error $\frac{\|\boldsymbol{M}_{\textrm{approx}} - \boldsymbol{M}_r\|_F}{\|\boldsymbol{M}_r \|_F}$.  }
 \end{subfigure}
 \begin{subfigure}[t]{0.5\textwidth}
	 \includegraphics[width=\textwidth]{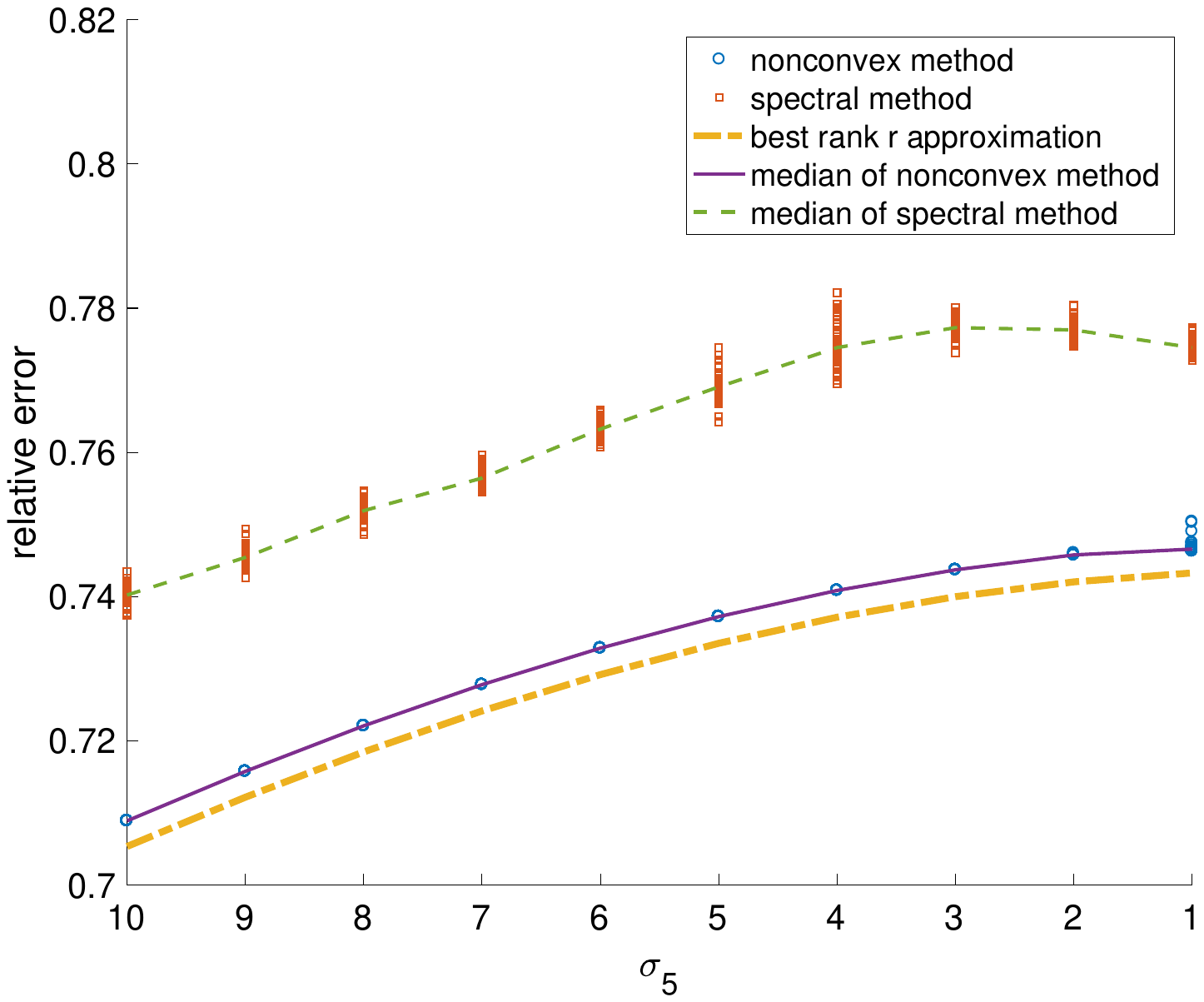}
	 \caption{    Relative error $\frac{\|\boldsymbol{M}_{\textrm{approx}} - \boldsymbol{M}\|_F}{\|\boldsymbol{M} \|_F}$.   }
 \end{subfigure}
 \caption{  Relative errors for full rank case.}
\label{small02_fg}
\end{figure}

Here $\boldsymbol{M}$ is assumed to have full rank, i.e., $\rank(\mtx{M}) = 500$. To be specific, let $\sigma_1 =\dots =  \sigma_4 = 10$, $\sigma_{6} = \dots = \sigma_{500} = 1$, and $\sigma_{5} = 10, 9, 8, \ldots, 2, 1$. The selected rank used in the nonconvex optimization \eqref{eq:obj_psd} is set as $r = 5$, and the sampling rate is set as $p=0.2$. With different values of $\sigma_{5}$, the results of our implementations of the gradient descent are plotted in Figure \ref{small02_fg}. One can observe that the relative errors for our nonconvex method \eqref{eq:obj_psd} are well-bounded for different $\sigma_5$'s, and much smaller than those for spectral low-rank approximation. The results indicate that our approach is able to approximate the ``true" best rank-$r$ approximation $\mtx{M}_r$ accurately in the presence of heavy spectral tail and possibly large condition number $\sigma_1/\sigma_{5}$, even with only $20\%$ observed entries.




\subsubsection{Low-rank matrix with large condition numbers}\label{small01_sim}
Here $\mtx{M}$ is assumed to be of exactly low rank with different condition numbers. Let $\sigma_1 =\dots =  \sigma_4 = 10$, $\sigma_{5} = \frac{10}{\kappa}$, and $\sigma_{6} = \dots = \sigma_{500} = 0$. Here the condition number takes on values $\kappa = 10,20,30,40,50,100,200,\infty$, which implies $\text{rank}(\mtx{M}) = 5$ if $\kappa < \infty$ while $\text{rank}(\mtx{M}) = 4$ if $\kappa = \infty$. The selected rank is always assumed to be $r = 5$, while the sampling rate is always $p=0.2$.

\begin{figure}
\center
\includegraphics[width = 0.5\textwidth]{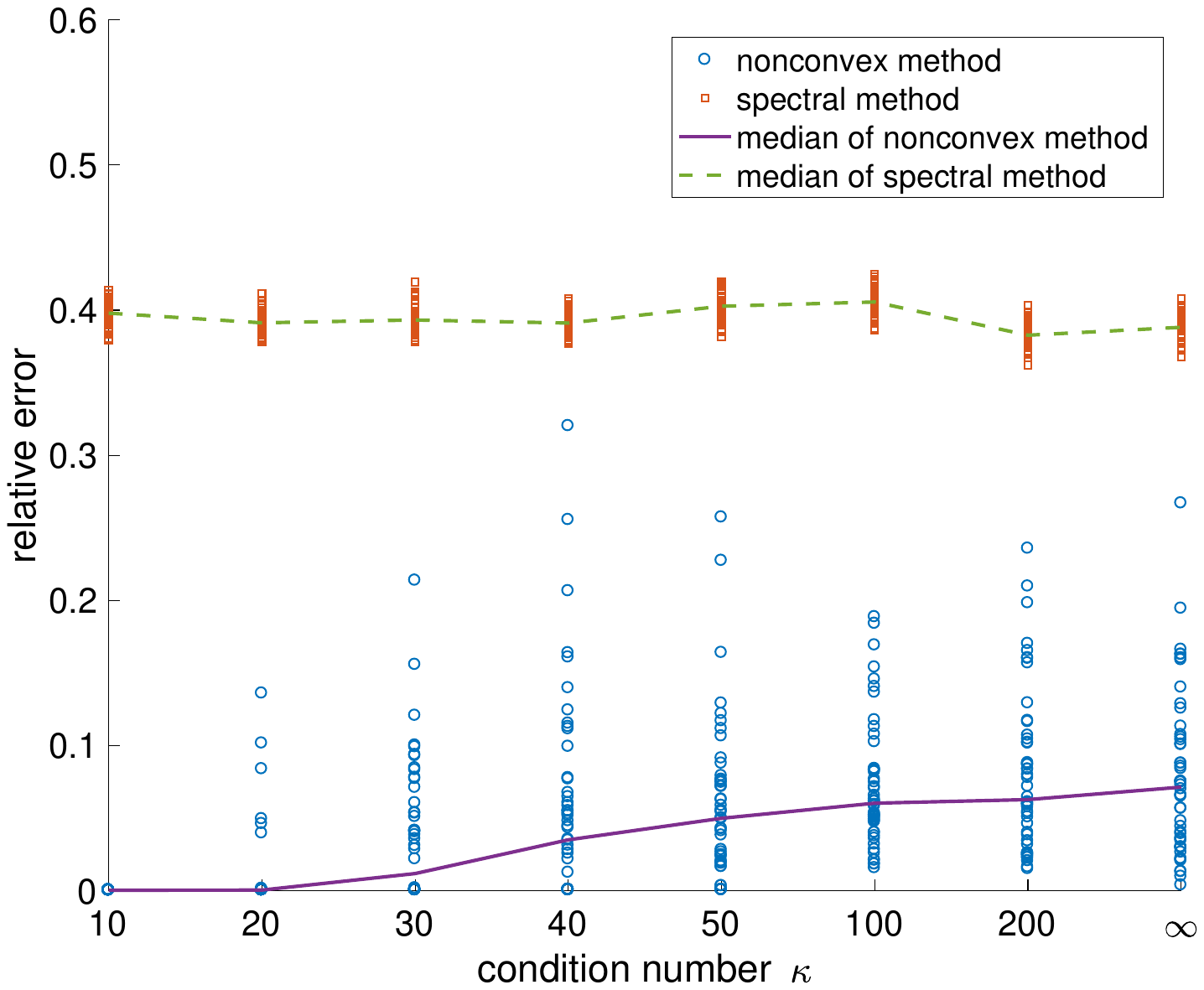}
\caption{Relative error $\frac{\|\boldsymbol{M}_{\textrm{approx}} - \boldsymbol{M}\|_F}{\|\boldsymbol{M} \|_F}$ for low-rank matrix with extreme condition numbers.}
\label{small01_fg}
\end{figure}

The performance of our nonconvex approach with various choices of $\kappa$ is demonstrated in Figure \ref{small01_fg}. One can observe that our nononvex optimization approach yields exact recovery of $\mtx{M}$ when $\kappa = 10$, while yields accurate low-rank approximation for $\mtx{M}$ with relative errors almost always smaller than $0.3$ when $\kappa \geqslant 20$. This fact is consistent with the example we discussed in Section \ref{sec:exmp1}, where we have shown that under certain incoherence conditions, the relative approximation error can be well-bounded even when $\kappa_r = \infty$.




\subsubsection{Rank mismatching}\label{small04_sim}
\begin{figure}[!ht]
 \begin{subfigure}[t]{0.5\textwidth}
	 \includegraphics[width=\textwidth]{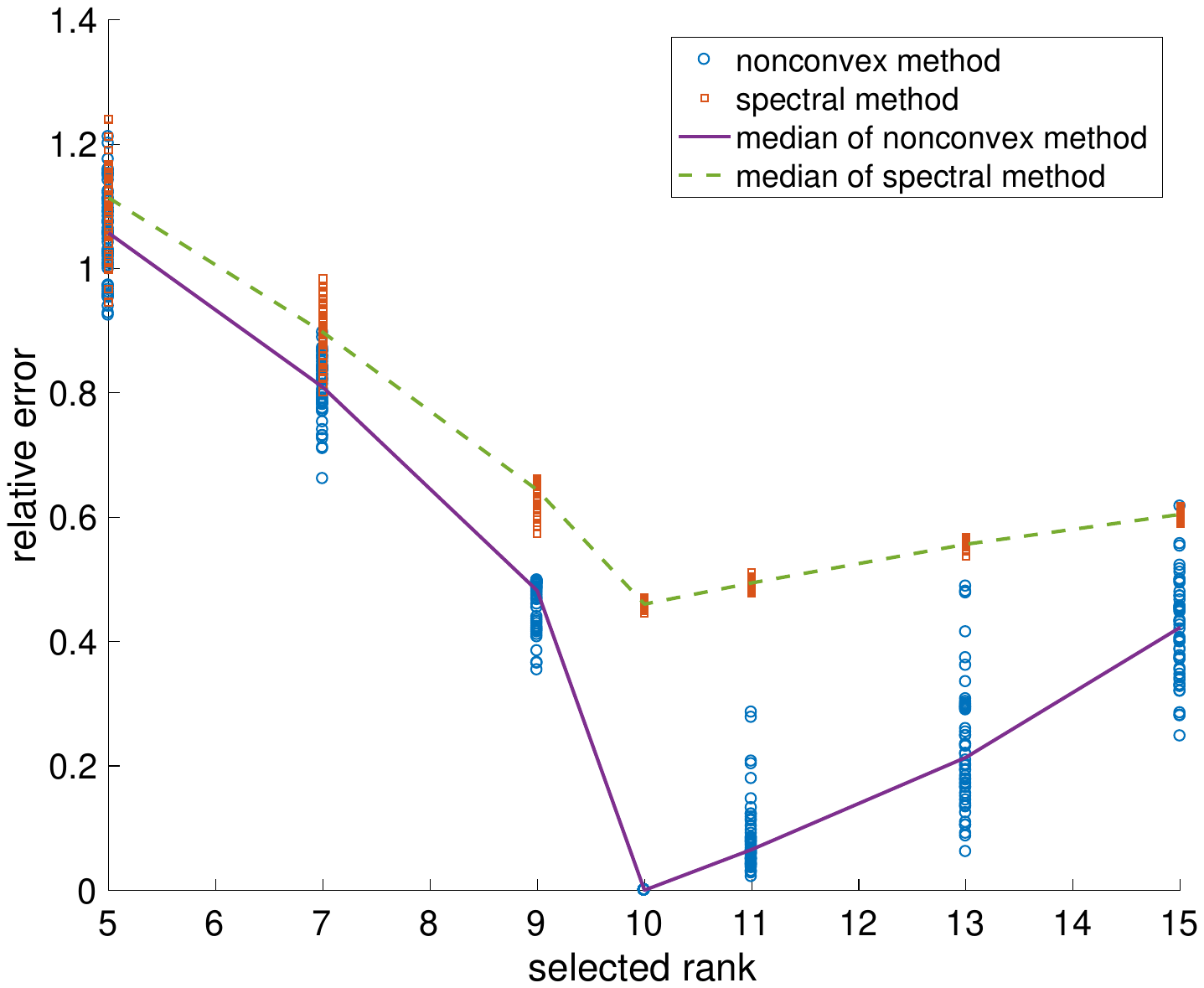}
	 \caption{   Relative error $\frac{\|\boldsymbol{M}_{\textrm{approx}} - \boldsymbol{M}_r\|_F}{\|\boldsymbol{M}_r \|_F}$.  }
 \end{subfigure}
 \begin{subfigure}[t]{0.5\textwidth}
	 \includegraphics[width=\textwidth]{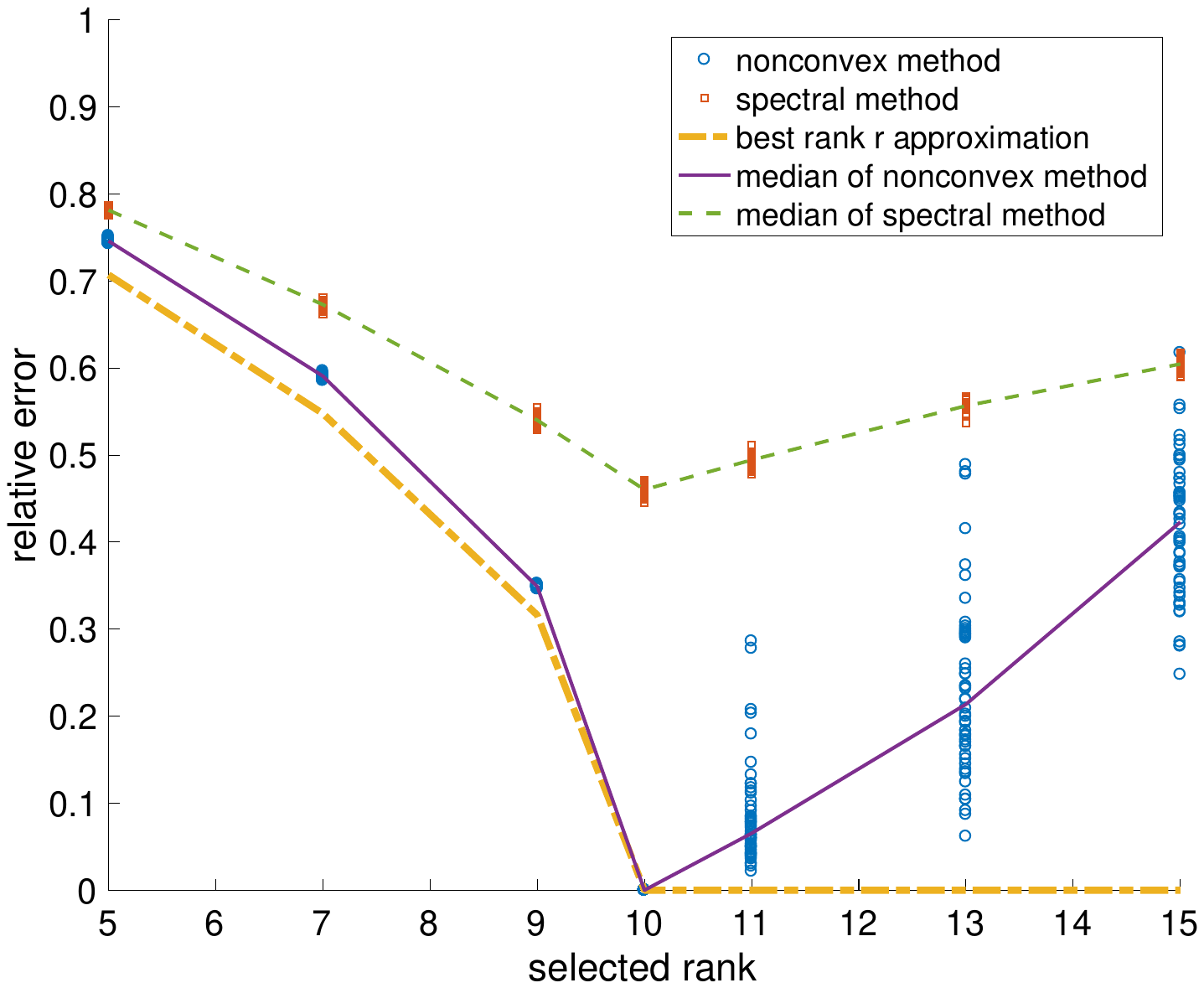}
	 \caption{   Relative error $\frac{\| \boldsymbol{M}_{\textrm{approx}} - \boldsymbol{M}\|_F}{\|\boldsymbol{M} \|_F}$.  }
 \end{subfigure}
 \caption{  Relative errors for rank mismatching for a fixed $\mtx{M}$ with $\text{rank}(\mtx{M})=10$.}
\label{small04_fg}
\end{figure}

\begin{figure}[!ht]
 \begin{subfigure}[t]{0.5\textwidth}
	 \includegraphics[width=\textwidth]{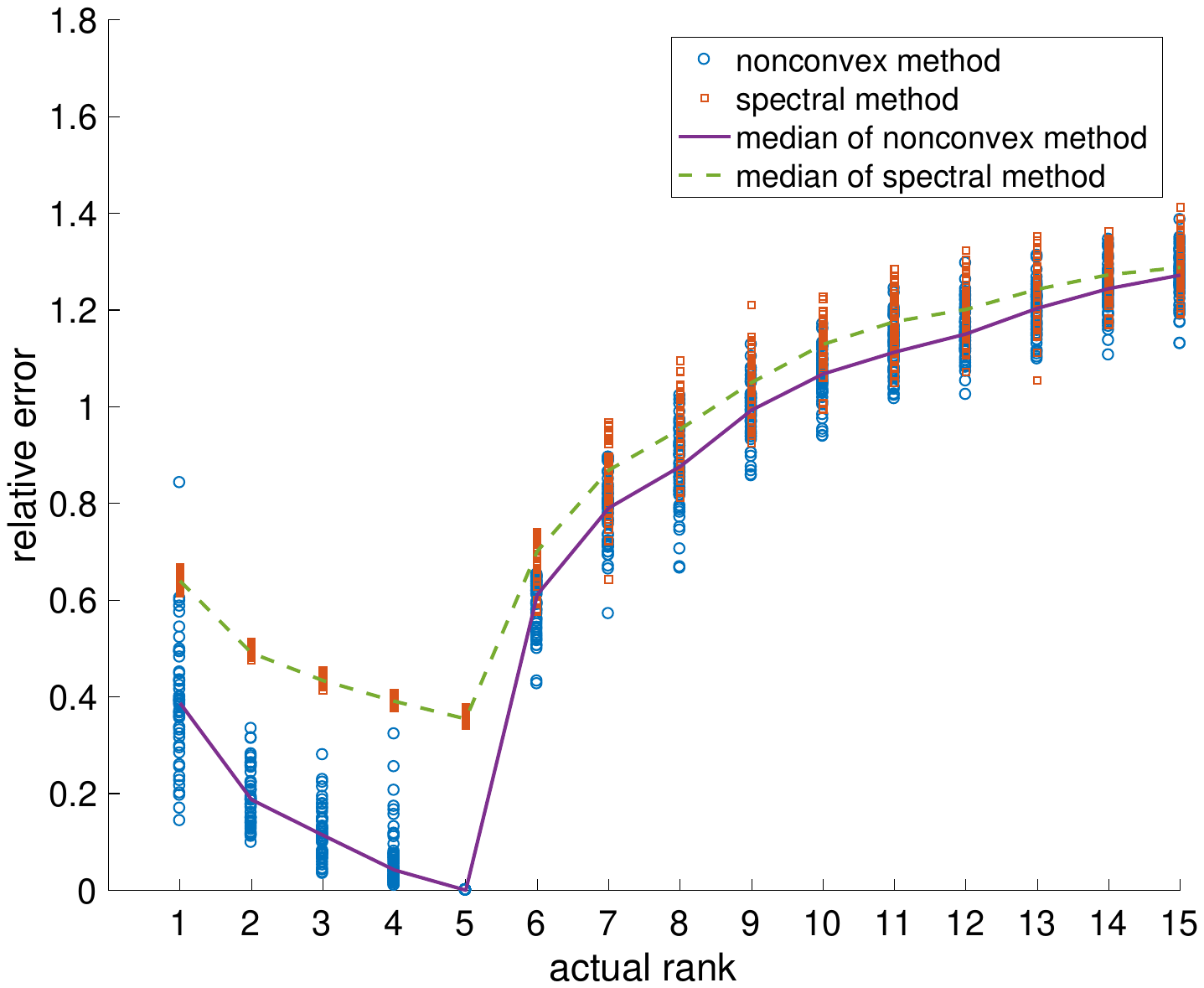}
	 \caption{   Relative error $\frac{\|\boldsymbol{M}_{\textrm{approx}} - \boldsymbol{M}_r\|_F}{\|\boldsymbol{M}_r \|_F}$.  }
 \end{subfigure}
 \begin{subfigure}[t]{0.5\textwidth}
	 \includegraphics[width=\textwidth]{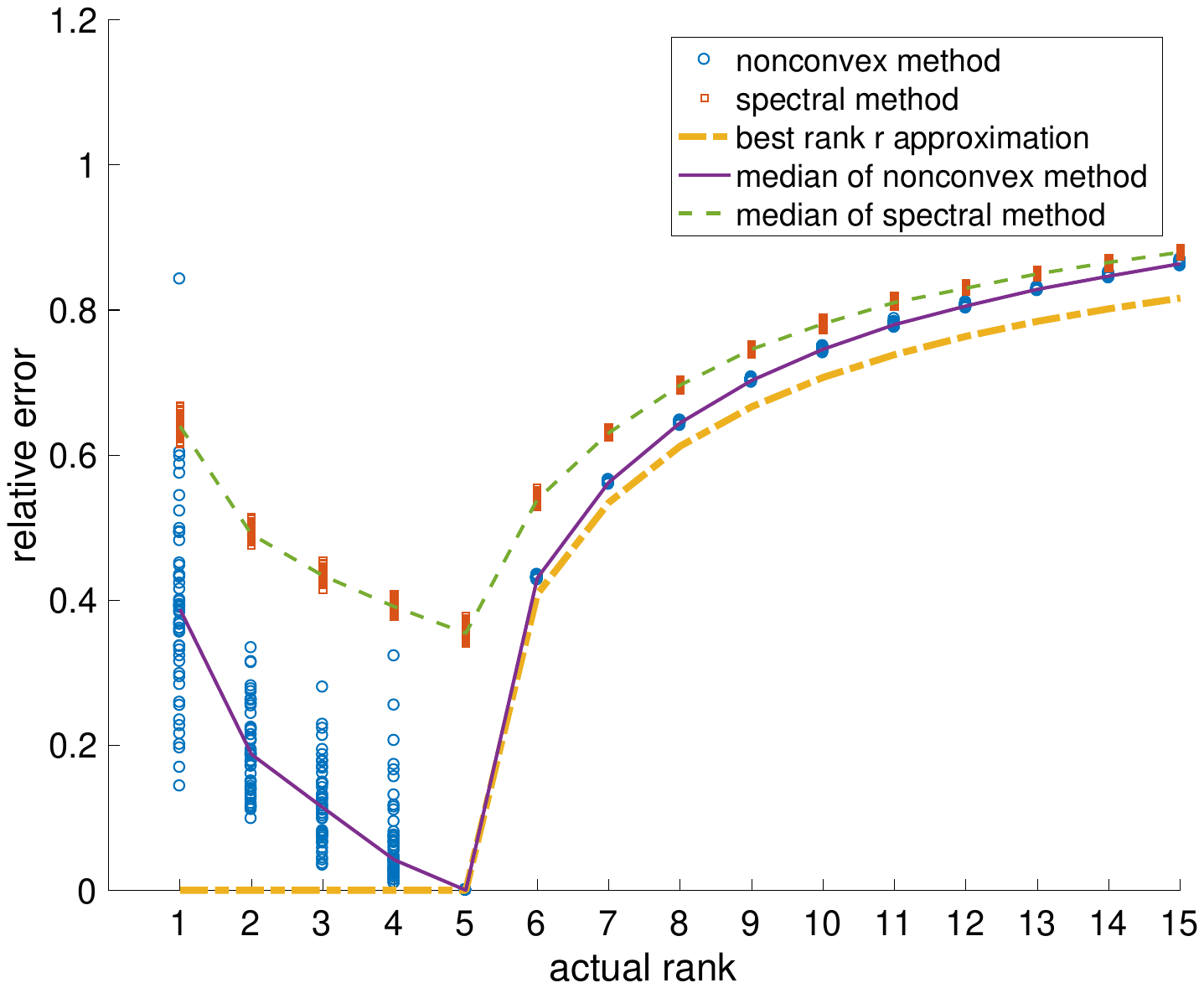}
	 \caption{   Relative error $\frac{\| \boldsymbol{M}_{\textrm{approx}} - \boldsymbol{M}\|_F}{\|\boldsymbol{M} \|_F}$.  }
 \end{subfigure}
 \caption{  Relative errors for rank mismatching, fixed selected rank.}
\label{small05_fg}
\end{figure}

\begin{figure}[!ht]
 \begin{subfigure}[t]{0.5\textwidth}
	 \includegraphics[width=\textwidth]{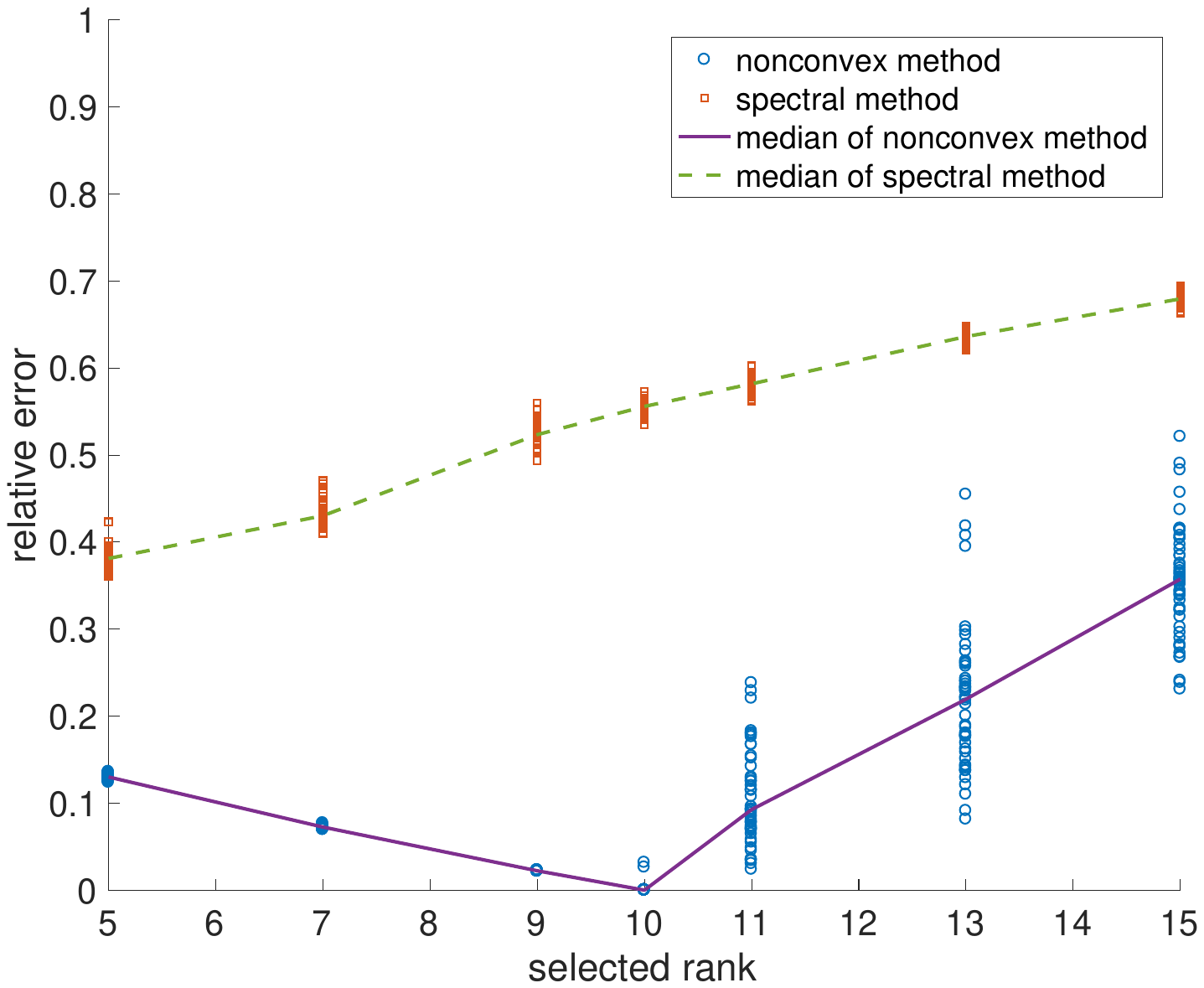}
	 \caption{   Relative error $\frac{\|\boldsymbol{M}_{\textrm{approx}} - \boldsymbol{M}_r\|_F}{\|\boldsymbol{M}_r \|_F}$.  }
 \end{subfigure}
 \begin{subfigure}[t]{0.5\textwidth}
	 \includegraphics[width=\textwidth]{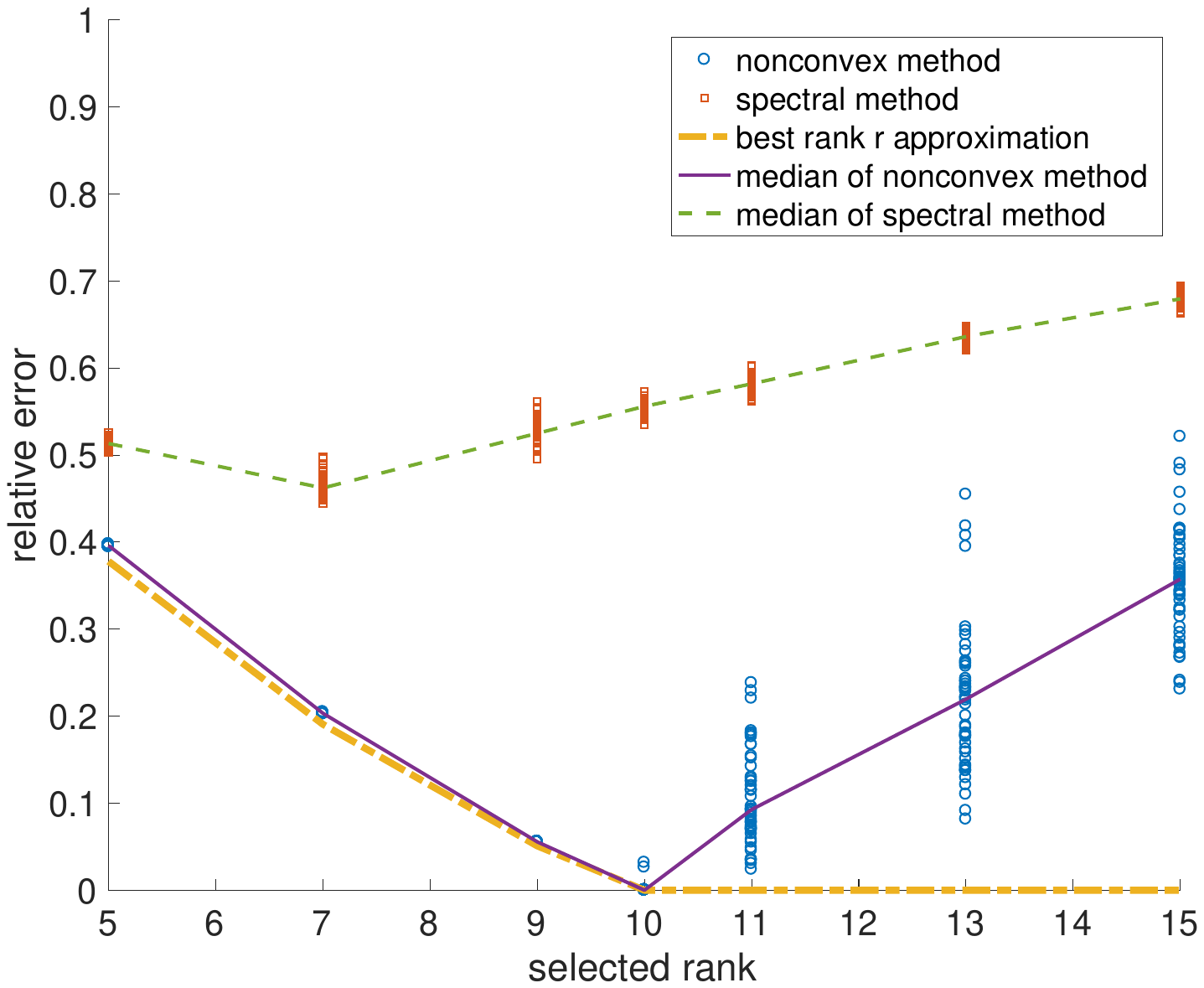}
	 \caption{   Relative error $\frac{\| \boldsymbol{M}_{\textrm{approx}} - \boldsymbol{M}\|_F}{\|\boldsymbol{M} \|_F}$.  }
 \end{subfigure}
 \caption{  Relative errors for rank mismatching for a fixed $\mtx{M}$ with $\text{rank}(\mtx{M})=10$.}
\label{neo_02_fg}
\end{figure}

\begin{figure}[!ht]
 \begin{subfigure}[t]{0.5\textwidth}
	 \includegraphics[width=\textwidth]{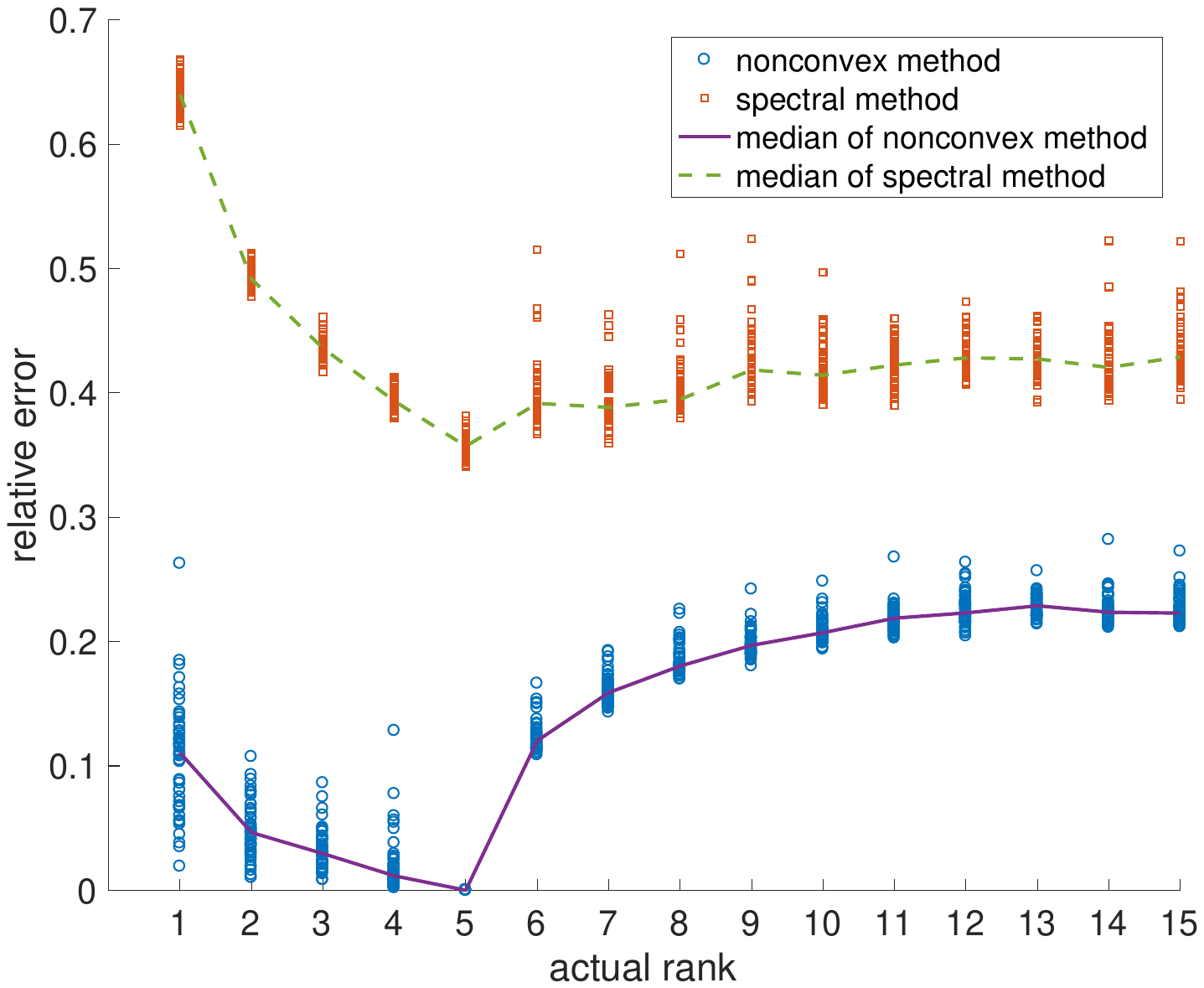}
	 \caption{   Relative error $\frac{\|\boldsymbol{M}_{\textrm{approx}} - \boldsymbol{M}_r\|_F}{\|\boldsymbol{M}_r \|_F}$.  }
 \end{subfigure}
 \begin{subfigure}[t]{0.5\textwidth}
	 \includegraphics[width=\textwidth]{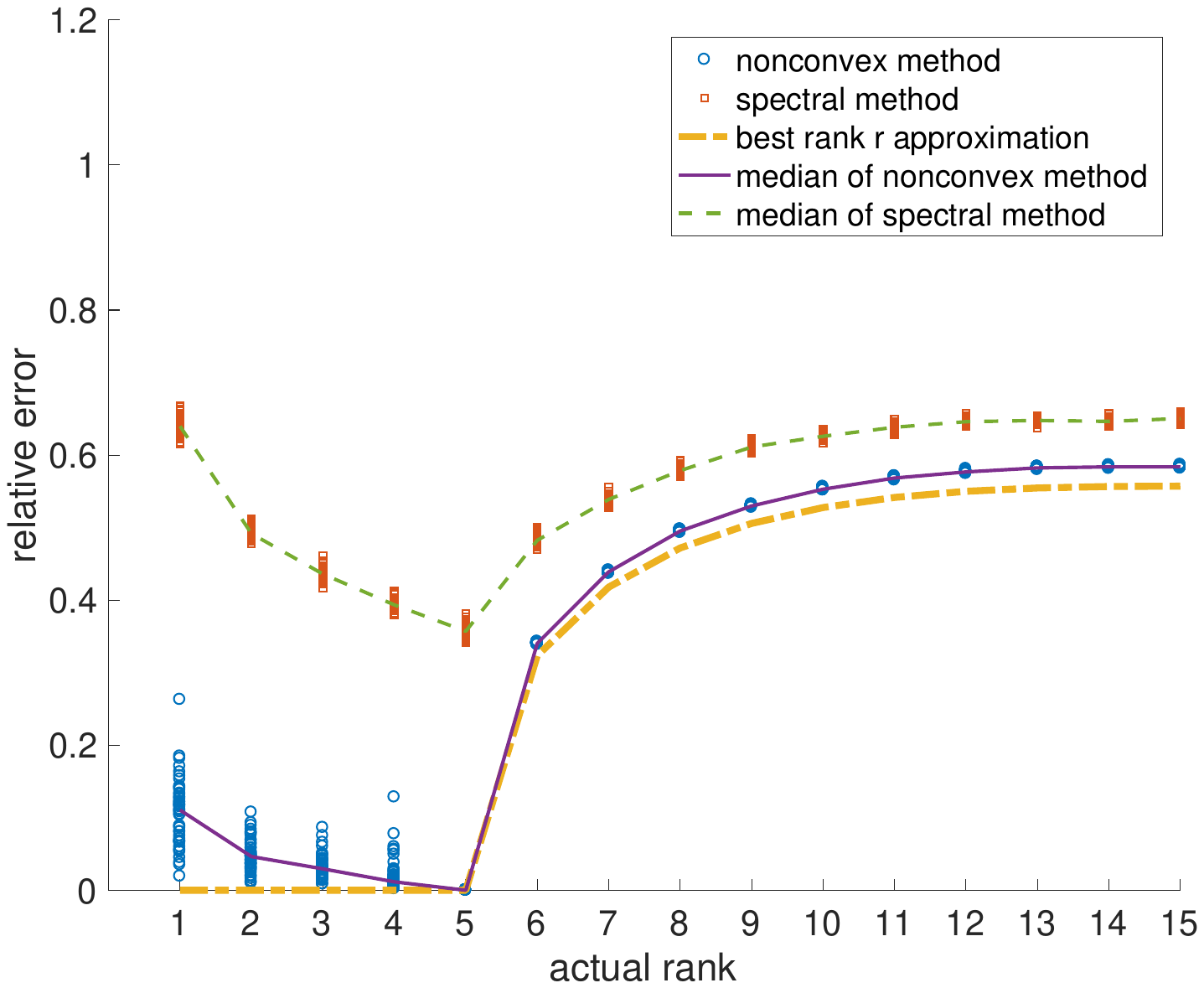}
	 \caption{   Relative error $\frac{\| \boldsymbol{M}_{\textrm{approx}} - \boldsymbol{M}\|_F}{\|\boldsymbol{M} \|_F}$.  }
 \end{subfigure}
 \caption{  Relative errors for rank mismatching, fixed selected rank.}
\label{neo_03_fg}
\end{figure}

In this section, we consider rank mismatching, i.e., the rank of $\mtx{M}$ is low but different from the selected rank $r$. In particular, we consider two settings for simulation: First, we fix $\boldsymbol{M}$ with $\text{rank}(\boldsymbol{M}) =10$, while the nonconvex optimization is implemented with selected rank $r = 5,7,9,10,11,13,15$; Second, the matrix $\boldsymbol{M}$ is randomly generated with rank from $1$ to $15$, while the selected rank is always $r=5$. The sampling rate is fixed as $p = 0.2$. We perform the simulation on two sets of spectrums: For the first one, all the nonzero eigenvalues are $10$; And the second one has decreasing eigenvalues: $\sigma_1 = 20,\sigma_2 = 18,\cdots,\sigma_{10} = 2$ for the case of fixed $\rank(\boldsymbol{M})$, $\sigma_1 = 30,\cdots,\sigma_{\rank(\boldsymbol{M})} = 32 - 2\times \rank(\boldsymbol{M})$ for the case of fixed selected rank $r$. Numerical results for the case of fixed $\rank(\mtx{M})$ are demonstrated in Figure \ref{small04_fg} (constant nonzero eigenvalues) and Figure \ref{neo_02_fg} (decreasing nonzero eigenvalues), while the case of fixed selected rank in Figure \ref{small05_fg} (constant nonzero eigenvalues) and Figure \ref{neo_03_fg} (decreasing nonzero eigenvalues). One can observe from these figures that if the selected rank $r$ is less than the actual rank $\text{rank}(\mtx{M})$, for the approximation of $\mtx{M}$, our nonconvex approach performs almost as well as the complete-data based best low-rank approximation $\mtx{M}_r$. Another interesting phenomenon is that our nonconvex method outperforms simple spectral methods in the approximation of either $\mtx{M}$ or $\mtx{M}_r$ significantly if the selected rank is greater than or equal to the true rank.

\subsection{Memory-efficient Kernel PCA}
\label{sec:KPCA_num}
In this section we study the empirical performance of our memory-efficient Kernel PCA approach by applying it to the synthetic data set in \citet{wang2012kernel}. The data set is an i.i.d. sample with sample size $n = 10,000$ and dimension $d=3$, and the data points are partitioned into two classes independently with equal probabilities. Points in the first class are first generated uniformly at random on the three-dimensional sphere $\{\vct{x}:\|\vct{x}\|_2=0.3\}$, while points in the second class are first generated uniformly at random on the three-dimensional sphere $\{\vct{x}:\|\vct{x}\|_2=1\}$. Every point is then perturbed independently by $\mathcal{N}(\boldsymbol{0}, \frac{1}{100}\boldsymbol{I}_3)$ noise. We aim to implement memory-efficient uncentered kernel PCA with $r=2$ on this dataset with the radial kernel $\exp(- \|\boldsymbol{x}-\boldsymbol{y}\|_2^2)$ in order to cluster the data points.

To implement the Nystr\"{o}m method \citep{williams2001using}, $50$ columns (and corresponding rows) are selected uniformly at random without replacement, then a rank-$2$ approximation of the kernel matrix $\boldsymbol{M}$ can be efficiently constructed with a smaller scale factorization. The effective sampling rate for Nystr\"{o}m method is $p_{\textrm{Nys}} = \frac{2\times 50n-50^2}{n^2}\approx 0.01$.
In contrast, in addition to recording the selected entry values, our nonconvex optimization method also requires to record the row and column indices for each selected entry. By using sparse matrix storage schemes like compressed sparse row (CSR) format \citep{saad2003iterative}, it needs $2n^2 p_{NCVX}+n+1$ entries to store the sparse matrix. Therefore, if $p_{NCVX} \geqslant \frac{3}{n}$, the nonconvex approach requires at most $2.5$ times as much memory as Nystr\"{o}m method for the same sampling complexity. Therefore, we choose the sampling rate $p_{\textrm{NCVX}} = \frac{p_{\textrm{Nys}}}{2.5}$ in the implementation of the nonconvex optimization \eqref{eq:obj_psd} such that the memory consumption is less costly than the Nystr\"{o}m method.

\begin{figure}[!ht] 
 \begin{subfigure}[t]{0.5\textwidth}
	 \includegraphics[width=\textwidth]{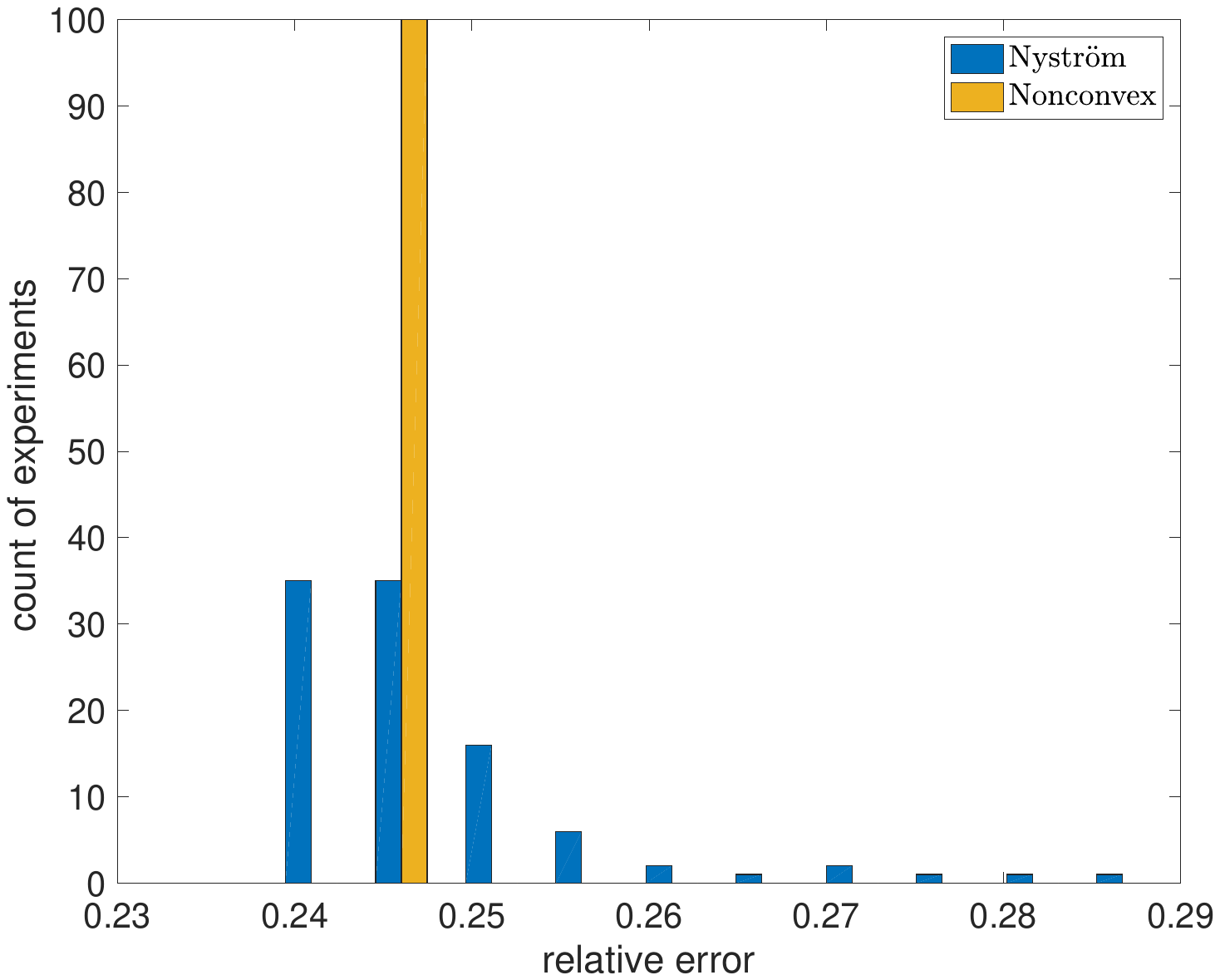}
\caption{  Relative error $\frac{\|\boldsymbol{M}_{\textrm{approx}} - \boldsymbol{M}_2\|_F}{\|\boldsymbol{M}_2 \|_F}$.  }
 \end{subfigure}
 \begin{subfigure}[t]{0.5\textwidth}
	 \includegraphics[width=\textwidth]{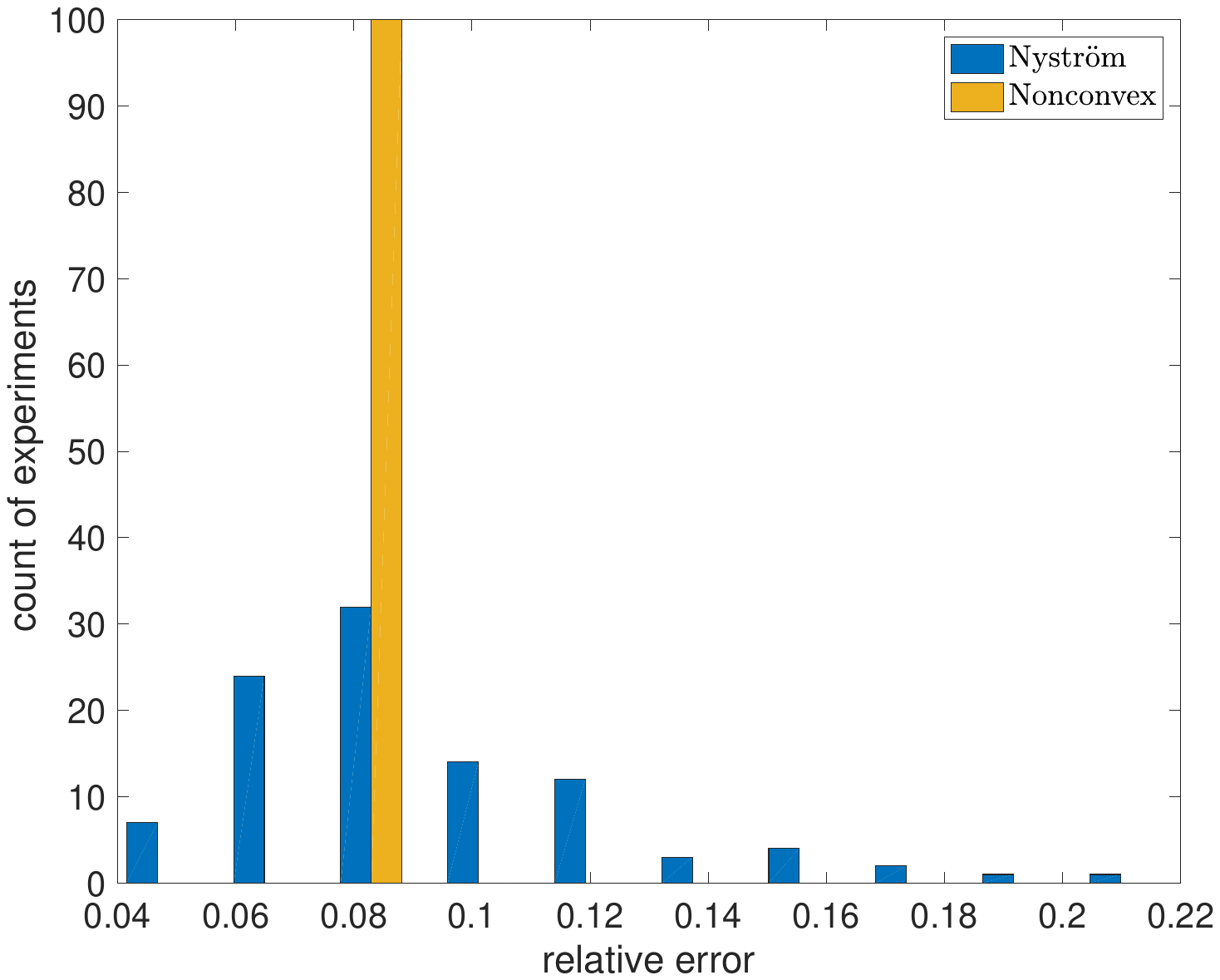}
 \caption{   Relative error $\frac{\|\boldsymbol{M}_{\textrm{approx}} - \boldsymbol{M}\|_F}{\|\boldsymbol{M} \|_F}$. }
 \end{subfigure}
\caption{  Relative errors for Nystr\"{o}m method with sampling rate $p_{\textrm{Nys}}\approx 0.01$ and nonconvex method with sampling rate $p_{\textrm{NCVX}}=\frac{p_{\textrm{Nys}}}{2.5}$.}
\label{large01_fg}
\end{figure}

Fixing such a synthetic data set, we apply both the Nystr\"{o}m method and our approach (with $\alpha = 100\|\boldsymbol{M}\|_{\ell_\infty} = 100$ and $\lambda = 500\sqrt{np_{\textrm{NCVX}}}$) for 100 times. Denote by $\boldsymbol{M}$ the ground truth of the kernel matrix, by $\boldsymbol{M}_2$ the ground truth of the best rank-$2$ approximation of $\boldsymbol{M}$, and by $\boldsymbol{M}_{\textrm{approx}}$ the memory efficient rank-$2$ approximation obtained by Nystr\"{o}m method or our nonconvex optimization. The left and right panels of Figure \ref{large01_fg} compare the two methods in approximating $\boldsymbol{M}_2$ and $\mtx{M}$ respectively based on the distributions of relative errors throughout the 100 Monte Carlo simulations. One can see that our approach is comparable with the Nystr\"{o}m method in terms of median performance, but much more stable. 

\begin{figure}[!ht] 
 \center
	 \includegraphics[width=0.5\textwidth]{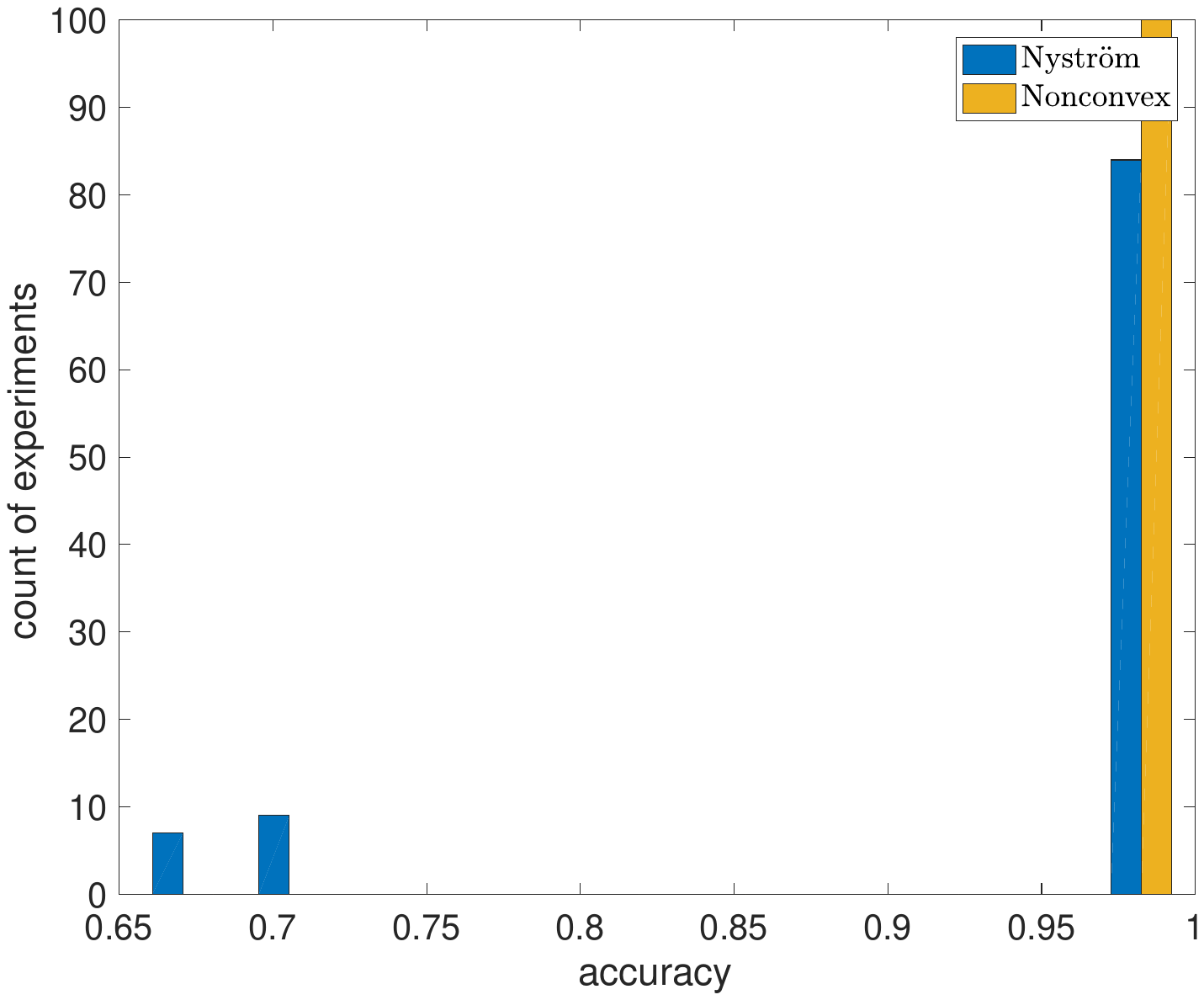}
\caption{  Clustering accuracy for Nystr\"{o}m method with sampling rate $p_{\textrm{Nys}}\approx 0.01$ and nonconvex method with sampling rate $p_{\textrm{NCVX}}=\frac{p_{\textrm{Nys}}}{2.5}$.  } 
\label{large02_fg}
\end{figure}

Both Nystr\"{o}m method and our nonconvex optimization \eqref{eq:obj_psd} give approximation in the form of $\mtx{M} \approx \widehat{\boldsymbol{X}}\widehat{\boldsymbol{X}}^\top$, so clustering analysis can be directly implemented based on $\widehat{\boldsymbol{X}}$. We implement k-means on the rows of $\widehat{\boldsymbol{X}}$ with $20$ repetitions, and Figure \ref{large02_fg} compares the two methods in the distribution of clustering accuracies. It clearly shows that our nonconvex optimization \eqref{eq:obj_psd} yields accurate clustering throughout the 100 tests while the Nystr\"{o}m method results in poor clustering occasionally. 

Moreover, during the iterations of the nonconvex method, the regularization term never activate throughout the 100 simulations. Therefore, empirically speaking, the performances of our numerical tests will remain the same if we simply set $\lambda = 0$.



\section{Proofs}
\label{sec:proofs}
In this section, we give a proof for main theorem. In Section \ref{sec_proofs_support}, we will present some useful supporting lemmas; in Section \ref{sec_proofs_main}, we present a proof for our main result Theorem  \ref{thm:main_psd}; finally in Section \ref{sec_proofs_proof} we give proof of lemmas used in former subsections. Our proof ideas benefit from those in \citet{ge2017no} as well as \citet{zhu2017global}, \citet{jin2017escape}.

\subsection{Supporting lemmas} 
\label{sec_proofs_support}

Here we give some useful supporting lemmas:

First, in literatures like \citet{vu2014simple} and \citet{bandeira2016sharp}, the control of $\|\boldsymbol{\Omega}-p\boldsymbol{J}\|$ is discussed:
\begin{lemma}\label{lemma_vu14_psd}
There is a constant $C>0$ such that the following holds. If $\Omega$ is sampled according to the off-diagonal symmetric $\text{Ber}(p)$ model, then
 \begin{equation*}
	 \mathbb{P}\left[\|\boldsymbol{\Omega}-p\boldsymbol{J}\|\geqslant  C \sqrt{np(1-p)} + C\sqrt{\log n}\right]\leqslant n^{-3}.
 \end{equation*}
\end{lemma}

~\\
Recall that eigen-space incoherence condition has been proposed in \citet{candes2009exact}.
\begin{definition}[\citealt{candes2009exact}]
 For any subspace $\mathcal{U}$ of $\mathbb{R}^n$ of dimension $r$, we define 
\begin{equation}
\label{eq:mu_subspace}
 \mu(\mathcal{U}) \coloneqq \frac{n}{r}\max_{1\leqslant i\leqslant n}\|\mathcal{P}_{\mathcal{U}}\boldsymbol{e}_i\|_2^2,
\end{equation}
where $\vct{e}_1, \ldots, \vct{e}_n$ represents the standard orthogonal basis of $\mathbb{R}^n$.	
\end{definition}

Similar to Theorem 4.1 in \citet{candes2009exact}, for the off-diagonal symmetric $\text{Ber}(p)$ model, we also have:

\begin{lemma}\label{lemma_rip_psd}
Suppose $\Omega$ is sampled according to the off-diagonal symmetric $\text{Ber}(p)$ model with probability $p$, define subspace
		 \begin{equation*}
			 \mathcal{T} \coloneqq \{\boldsymbol{M}\in \mathbb{R}^{n \times n}\mid (\mathcal{I}-\mathcal{P}_{\mathcal{U}})\boldsymbol{M}(\mathcal{I}-\mathcal{P}_{\mathcal{U}}) = \boldsymbol{0},\;\boldsymbol{M}\;\textrm{symmetric}\},
		 \end{equation*}
		 where $\mathcal{U}$ is a fixed subspace of $\mathbb{R}^n$. Let $\mathcal{P}_{\mathcal{T}}$ be the Euclidean projection on to $\mathcal{T}$. Then there is an absolute constant $C$, for any $\delta\in(0,1]$, if $p\geqslant  C\frac{\mu(\mathcal{U}) \dim(\mathcal{U})\log n}{\delta^2 n}$ with $\mu(\mathcal{U})$ defined in \eqref{eq:mu_subspace}, in an event $E$ with probability $\mathbb{P}[E]\geqslant  1-n^{-3}$, we have
		 \begin{equation*}
			 p^{-1}\|\mathcal{P}_{\mathcal{T}} \mathcal{P}_\Omega \mathcal{P}_{\mathcal{T}}-p \mathcal{P}_{\mathcal{T}}\|\leqslant \delta .
		 \end{equation*}  
\end{lemma}
In \citet{gross2011recovering} and \citet{gross2010note}, similar results are given for symmetric uniform sampling with/without replacement. The proof of Lemma \ref{lemma_rip_psd} are very similar to those in \citet{recht2011simpler}.

~\\
For the first and second order optimality condition of our objective function $f(\boldsymbol{X})$ defined in \eqref{eq:obj_psd}, we have:
\begin{lemma}[{\citealt[Proposition 4.1]{ge2016matrix}}]\label{lemma_12condition}
The first order optimality condition of objective function \eqref{eq:obj_psd} is
\begin{equation*}
 \nabla f(\boldsymbol{X}) = 2\mathcal{P}_{\Omega}(\boldsymbol{X}\boldsymbol{X}^\top-\boldsymbol{M})\boldsymbol{X} + \lambda \nabla G_{\alpha}(\boldsymbol{X}) = \boldsymbol{0},
\end{equation*}
and the second order optimality condition requires $\forall \boldsymbol{H}\in\mathbb{R}^{n\times r}$, we have
\begin{equation*}
\begin{split}
&\operatorname{vec}(\boldsymbol{H})^\top \nabla^2 f(\boldsymbol{X})\operatorname{vec}(\boldsymbol{H}) \\
 =& \|\mathcal{P}_{\Omega}(\boldsymbol{H}\boldsymbol{X}^\top+\boldsymbol{X}\boldsymbol{H}^\top)\|_F^2 + 2\langle \mathcal{P}_{\Omega}(\boldsymbol{X}\boldsymbol{X}^\top - \boldsymbol{M}),\mathcal{P}_{\Omega}(\boldsymbol{H}\boldsymbol{H}^\top) \rangle +  {\lambda}\operatorname{vec}(\boldsymbol{H})^\top \nabla^2 G_{\alpha}(\boldsymbol{X})\operatorname{vec}(\boldsymbol{H})\\
 \geqslant & 0.
\end{split}
\end{equation*}
\end{lemma}

~\\
Finally, we are going to present our main lemma which will be used multiple times. Before we formally state it, for simplicity of notations, for any matrix $\boldsymbol{M}_1,\boldsymbol{M}_2\in\mathbb{R}^{n_1\times n_2}$, any set $\Omega_0 \in[n_1]\times [n_2]$ and any real number $t\in\mathbb{R}$, we introduce following notation:
\begin{equation}\label{eq:DOmegap}
D_{\Omega_0,t}(\boldsymbol{M}_1,\boldsymbol{M}_2) \coloneqq \langle \mathcal{P}_{\Omega_0}(\boldsymbol{M}_1),\mathcal{P}_{\Omega_0}(\boldsymbol{M}_2)\rangle - t\langle \boldsymbol{M}_1, \boldsymbol{M}_2\rangle.
\end{equation}
Now we are well prepared to present our main lemma:

\begin{lemma}
\label{lemma_concern} 
 For any $\Omega_0\subset [n_1]\times [n_2]$ and corresponding $\boldsymbol{\Omega}_0$, for all $\boldsymbol{A}\in\mathbb{R}^{n_1\times r_1},\boldsymbol{B}\in \mathbb{R}^{n_1\times r_2},\boldsymbol{C}\in\mathbb{R}^{n_2\times r_1},\boldsymbol{D}\in \mathbb{R}^{n_2\times r_2}$, we have
 \begin{equation}\label{eq_concern}
 |D_{\Omega_0,t}(\boldsymbol{A}\boldsymbol{C}^\top, \boldsymbol{B}\boldsymbol{D}^\top)| \leqslant  \|\boldsymbol{\Omega}_0-t\boldsymbol{J}\| \sqrt{\sum_{k=1}^{n_1}\|\boldsymbol{A}_{k,\cdot}\|_2^2 \|\boldsymbol{B}_{k,\cdot}\|_2^2}\sqrt{\sum_{k=1}^{n_2}\|\boldsymbol{C}_{k,\cdot}\|_2^2 \|\boldsymbol{D}_{k,\cdot}\|_2^2}
 \end{equation}
 for all $t\in\mathbb{R}$.
\end{lemma}

We will use this result for $\Omega_0 = \Omega,t = p$ multiple times later. Note here we do not have any assumption on $\Omega_0$ and this is a deterministic result. The proof of this lemma can be considered as a more sophisticated version of those proofs for spectral lemmas in \citet{bhojanapalli2014universal} and \citet{li2016recovery}:

\begin{lemma}[\citealt{bhojanapalli2014universal,li2016recovery}]\label{lemma_spec}
 Suppose matrix $\boldsymbol{M}\in\mathbb{R}^{n_1\times n_2}$ can be decomposed as $\boldsymbol{M} = \boldsymbol{B}\boldsymbol{D}^\top$, let $\Omega\subset [n_1]\times [n_2]$ be set of revealed entries (not necessary follow any specific distribution), then for any $t$ we have
 \begin{equation*}
	 \|\mathcal{P}_{\Omega}(\boldsymbol{M})-t\boldsymbol{M}\|\leqslant \|\boldsymbol{\Omega}-t\boldsymbol{J}\|\|\boldsymbol{B}\|_{2,\infty}\|\boldsymbol{D}\|_{2,\infty}.
 \end{equation*}
\end{lemma}

~\\
Lemma \ref{lemma_concern} is applied in our proof in replace of Theorem D.1 in \citet{ge2016matrix} to derive tighter inequalities. For comparison, we give the statement of that result:

\begin{lemma}[{\citealt[Theorem D.1]{ge2016matrix}}]\label{GLM16_thm_d.1}
 With high probability over the choice of $\Omega$, for any two rank-$r$ matrices $\boldsymbol{W},\boldsymbol{Z}\in\mathbb{R}^{n\times n}$, we have
 \begin{equation*}
	 \begin{split}
	 &|\langle \mathcal{P}_{\Omega}(\boldsymbol{W}),\mathcal{P}_{\Omega}(\boldsymbol{Z})  \rangle - p\langle \boldsymbol{W},\boldsymbol{Z} \rangle|\\
	 \leqslant& O\left(\|\boldsymbol{W}\|_{\ell_\infty}\|\boldsymbol{Z}\|_{\ell_\infty} nr\log n+\sqrt{pnr \|\boldsymbol{W}\|_{\ell_\infty}\|\boldsymbol{Z}\|_{\ell_\infty} \|\boldsymbol{W}\|_{F}\|\boldsymbol{Z}\|_{F}\log n}\right).
	 \end{split}
 \end{equation*}
\end{lemma}

~\\
In \citet{sun2015guaranteed}, \citet{chen2015fast} and \citet{zheng2016convergence}, the authors give an upper bound of $\|\mathcal{P}_{\Omega}(\boldsymbol{H}\boldsymbol{H}^\top)\|_F^2$ for any $\boldsymbol{H}$. To be more precise, they assume $\Omega$ is sampled according to i.i.d. Bernoulli model with probability $p$, then if $p\geqslant  C_1\frac{\log n}{n}$ with absolute constant $C_1$ sufficient large, 
 \begin{equation}\label{eq_015}
	 \|\mathcal{P}_{\Omega}(\boldsymbol{H}\boldsymbol{H}^\top)\|_F^2 \leqslant p\|\boldsymbol{H}\|_F^4 + C_2\sqrt{np} \sum_{i=1}^n \|\boldsymbol{H}_{i,\cdot}\|_2^4
 \end{equation}
 holds with high probability. 
 
 In contrast, by applying Lemma \ref{lemma_vu14_psd} and Lemma \ref{lemma_concern},
 \begin{equation}\label{eq_016}
	 |\|\mathcal{P}_{\Omega}(\boldsymbol{H}\boldsymbol{H}^\top)\|_F^2 - p\|\boldsymbol{H}\boldsymbol{H}^\top\|_F^2| \leqslant C_3\sqrt{np} \sum_{i=1}^n \|\boldsymbol{H}_{i,\cdot}\|_2^4
 \end{equation}
 holds with high probability, which is a tighter bound once we notice the fact that $\|\boldsymbol{H}\boldsymbol{H}^\top\|_F\leqslant \|\boldsymbol{H}\|_F^2 $. Moreover, comparing to \eqref{eq_015}, our result \eqref{eq_016} can measure the difference between $\|\mathcal{P}_{\Omega}(\boldsymbol{H}\boldsymbol{H}^\top)\|_F^2$ and its expectation $p\|\boldsymbol{H}\boldsymbol{H}^\top\|_F^2$, which makes the model-free analysis possible.
 

\subsection{A proof of Theorem \ref{thm:main_psd}}
\label{sec_proofs_main}
Here we give a proof of Theorem \ref{thm:main_psd}. The proof will be mainly divided into two parts: In Section \ref{sec_proofs_main_step1}, we discuss the landscape of objective function $f(\boldsymbol{X})$ and then define the auxiliary function $K(\boldsymbol{X})$, one can see the span of local minima of $f(\boldsymbol{X})$ can be controlled by the span of supperlevel set of $K(\boldsymbol{X})$: $\{\boldsymbol{X}\in\mathbb{R}^{n\times r}\mid K(\boldsymbol{X})\geqslant  0\}$; In Section \ref{sec_proofs_main_step3}, we give a uniform upper bound of $K(\boldsymbol{X})$ and use it to solve for possible span of supperlevel set of $K(\boldsymbol{X})$, which finishes the proof of our main result.

\subsubsection{Landscape of objective function $f$ and the auxiliary function $K$}
\label{sec_proofs_main_step1}
Now denote $\boldsymbol{U}_r \coloneqq [\sqrt{\sigma_1}\boldsymbol{u}_1\; \dots\; \sqrt{\sigma_r}\boldsymbol{u}_r]$. For a given $\boldsymbol{X}\in\mathbb{R}^{n\times r}$, suppose $\boldsymbol{X}^\top\boldsymbol{U}_r $ has SVD $\boldsymbol{X}^\top\boldsymbol{U}_r  = \boldsymbol{A}\boldsymbol{D}\boldsymbol{B}^\top$, and let $\boldsymbol{R}_{\boldsymbol{X},\boldsymbol{U}_r} \coloneqq \boldsymbol{B}\boldsymbol{A}^\top\in\mathsf{O}(r)$ and $\boldsymbol{U} \coloneqq \boldsymbol{U}_r\boldsymbol{R}_{\boldsymbol{X},\boldsymbol{U}_r}$, where $\mathsf{O}(r)$ denotes the set of $r\times r$ orthogonal matrices $\{\boldsymbol{R}\in\mathbb{R}^{r\times r}\mid\boldsymbol{R}^\top\boldsymbol{R} = \boldsymbol{R}\boldsymbol{R}^\top = \boldsymbol{I}\}$. Then one can verify that $\boldsymbol{X}^\top\boldsymbol{U} = \boldsymbol{A}\boldsymbol{D}\boldsymbol{A}^\top$ is a positive semidefinite matrix. Notice by the way we define $\boldsymbol{U}_r$ and $\boldsymbol{U}$, $\boldsymbol{U}_r\boldsymbol{U}_r^\top = \boldsymbol{U}\boldsymbol{U}^\top$. Moreover, $\boldsymbol{R}_{\boldsymbol{X},\boldsymbol{U}_r} \in \argmin_{\boldsymbol{R}\in\mathsf{O}(r)}\|\boldsymbol{X} - \boldsymbol{U}_r\boldsymbol{R}\|_F$, see, e.g., \citet{chen2015fast}.

~\\
Let $\boldsymbol{\Delta} \coloneqq \boldsymbol{X}-\boldsymbol{U}$, and define the following auxiliary function introduced in \citet{jin2017escape} and \citet{ge2017no}:
\begin{equation*}
 K(\boldsymbol{X}) \coloneqq \operatorname{vec}(\boldsymbol{\Delta})^\top \nabla^2 f(\boldsymbol{X})\operatorname{vec}(\boldsymbol{\Delta})-4 \langle \nabla f(\boldsymbol{X}),\boldsymbol{\Delta}\rangle.
\end{equation*}
The first and second order optimality conditions for any local minimum $\widehat{\boldsymbol{X}}$ imply that $K(\widehat{\boldsymbol{X}} ) \geqslant  0$. In other words, we have
\begin{equation*}
 \{\text{All local minima of~}  f(\boldsymbol{X})\}\subset \{\boldsymbol{X}\in \mathbb{R}^{n\times r} \mid K(\boldsymbol{X})\geqslant  0\}.
\end{equation*}
To study the properties of the local minima of $f(\boldsymbol{X})$, we can consider the supperlevel set of $K(\boldsymbol{X})$: $\{\boldsymbol{X}\in \mathbb{R}^{n\times r} \mid K(\boldsymbol{X})\geqslant  0\}$ instead. In order to get a clear representation of $K(\boldsymbol{X})$, one can plug in first and second order condition listed in Lemma \ref{lemma_12condition}. Actually, by repacking terms in \citet[Lemma 7]{ge2017no} and noticing the fact that by definition we have $\langle \boldsymbol{U}\boldsymbol{\Delta}^\top,\boldsymbol{N}\rangle = 0$, we can decompose $K(\boldsymbol{X})$ as following:
\begin{lemma}[{\citealt[Lemma 7]{ge2017no}}]\label{lemma_k_psd}
Uniformly for all $\boldsymbol{X}\in\mathbb{R}^{n\times r}$ and corresponding $\boldsymbol{\Delta}$ defined as $\boldsymbol{\Delta} \coloneqq \boldsymbol{X}-\boldsymbol{U}$, we have

\begin{equation}\label{eq:K}
 \begin{split}
 K(\boldsymbol{X}) =& \underbrace{p\left(\|\boldsymbol{\Delta}\boldsymbol{\Delta}^\top\|_F^2-3\|\boldsymbol{X}\boldsymbol{X}^\top -\boldsymbol{U}\boldsymbol{U}^\top\|_F^2 \right)}_{K_1(\boldsymbol{X})} \\
 &+ \underbrace{ D_{\Omega,p}(\boldsymbol{\Delta}\boldsymbol{\Delta}^\top,\boldsymbol{\Delta}\boldsymbol{\Delta}^\top) - 3D_{\Omega,p}(\boldsymbol{X}\boldsymbol{X}^\top-\boldsymbol{U}\boldsymbol{U}^\top,\boldsymbol{X}\boldsymbol{X}^\top-\boldsymbol{U}\boldsymbol{U}^\top)  }_{K_2(\boldsymbol{X})}\\
	 &+\underbrace{\lambda \left( \operatorname{vec}(\boldsymbol{\Delta})^\top\nabla^2 G_{\alpha}(\boldsymbol{X})\operatorname{vec}(\boldsymbol{\Delta})-4\langle \nabla G_{\alpha}(\boldsymbol{X}),\boldsymbol{\Delta}\rangle\right)}_{K_3(\boldsymbol{X})} \\
	 &+\underbrace {6D_{\Omega,p}( \boldsymbol{\Delta}\boldsymbol{\Delta}^\top,\boldsymbol{N} )+8D_{\Omega,p}( \boldsymbol{U}\boldsymbol{\Delta}^\top, \boldsymbol{N} ) + 6p \langle \boldsymbol{\Delta}\boldsymbol{\Delta}^\top,\boldsymbol{N} \rangle}_{K_4(\boldsymbol{X})},
\end{split}
\end{equation}
where $D_{\Omega,p}(\cdot,\cdot)$ is defined in \eqref{eq:DOmegap}.
\end{lemma}

Notice in Theorem \ref{thm:main_psd}, we are only concerned about the difference between $\boldsymbol{X}\boldsymbol{X}^\top$ and $\boldsymbol{M}_r$ (or $\boldsymbol{M}$), thus there is no difference to consider $\boldsymbol{X}$ or $\widetilde{\boldsymbol{X}} = \boldsymbol{X}\boldsymbol{R},\forall \boldsymbol{R}\in\mathsf{O}(r)$. Moreover, by the definition of $K(\boldsymbol{X})$, we have $K(\boldsymbol{X}) = K(\widetilde{\boldsymbol{X}})$. 

In fact, by the definition of $\boldsymbol{R}_{\boldsymbol{X},\boldsymbol{U}_r}$, we have $\boldsymbol{R}_{\boldsymbol{X}\boldsymbol{R},\boldsymbol{U}_r} = \boldsymbol{R}_{\boldsymbol{X},\boldsymbol{U}_r}\boldsymbol{R},\forall \boldsymbol{R}\in\mathsf{O}(r)$, which implies $\widetilde{\boldsymbol{U}} = \boldsymbol{U}\boldsymbol{R}$ and $\widetilde{\boldsymbol{\Delta}} = \boldsymbol{\Delta}\boldsymbol{R}$. Now we have
$$
 \widetilde{\boldsymbol{X}}\widetilde{\boldsymbol{X}}^\top = \boldsymbol{X}\boldsymbol{X}^\top, \widetilde{\boldsymbol{U}}\widetilde{\boldsymbol{U}}^\top = \boldsymbol{U}\boldsymbol{U}^\top, \widetilde{\boldsymbol{\Delta}}\widetilde{\boldsymbol{\Delta}}^\top = \boldsymbol{\Delta}\boldsymbol{\Delta}^\top,\widetilde{\boldsymbol{U}}\widetilde{\boldsymbol{\Delta}}^\top = \boldsymbol{U}\boldsymbol{\Delta}^\top,
$$
which means $K_i(\widetilde{\boldsymbol{X}}) = K_i( \boldsymbol{X})$ for $i = 1,2,4$. As for $K_3$, by \citet[Lemma 18]{ge2017no}, we have 
\begin{equation*}
 \begin{split}
	 &\operatorname{vec}(\boldsymbol{\Delta})^\top \nabla^2 G_{\alpha}(\boldsymbol{X})\operatorname{vec}(\boldsymbol{\Delta}) - 4\langle \nabla G_{\alpha}(\boldsymbol{X}),\boldsymbol{\Delta}\rangle\\
		=& 4\sum_{i=1}^n [(\|\boldsymbol{X}_{i,\cdot}\|_2-\alpha)_+]^3 \frac{\|\boldsymbol{X}_{i,\cdot}\|_2^2\|\boldsymbol{\Delta}_{i,\cdot}\|_2^2 -\langle \boldsymbol{X}_{i,\cdot},\boldsymbol{\Delta}_{i,\cdot}\rangle^2}{\|\boldsymbol{X}_{i,\cdot}\|_2^3}+ 12\sum_{i=1}^n [(\|\boldsymbol{X}_{i,\cdot}\|_2-\alpha)_+]^2 \frac{\langle \boldsymbol{X}_{i,\cdot},\boldsymbol{\Delta}_{i,\cdot}\rangle^2}{\|\boldsymbol{X}_{i,\cdot}\|_2^2}\\
		&   -16\sum_{i=1}^n [(\|\boldsymbol{X}_{i,\cdot}\|_2-\alpha)_+]^3 \frac{\langle \boldsymbol{X}_{i,\cdot},\boldsymbol{\Delta}_{i,\cdot}\rangle}{\|\boldsymbol{X}_{i,\cdot}\|_2}.
 \end{split}
\end{equation*}
Since $\boldsymbol{R}\in\mathsf{O}(r)$, we have $\|\widetilde{\boldsymbol{X}}_{i,\cdot}\|_2 = \|\boldsymbol{X}_{i,\cdot}\|_2$, $\|\widetilde{\boldsymbol{\Delta}}_{i,\cdot}\|_2 = \|\boldsymbol{\Delta}_{i,\cdot}\|_2$ and $\langle \widetilde{\boldsymbol{X}}_{i,\cdot},\widetilde{\boldsymbol{\Delta}}_{i,\cdot} \rangle = \langle\boldsymbol{X}_{i,\cdot},\boldsymbol{\Delta}_{i,\cdot} \rangle$, so we have $K_3(\widetilde{\boldsymbol{X}}) = K_3( \boldsymbol{X})$. Putting things together, we have $K(\widetilde{\boldsymbol{X}}) = K( \boldsymbol{X})$. 

Therefore, if we want to show any $\boldsymbol{X}$ with $K(\boldsymbol{X})\geqslant  0$ satisfies \eqref{eq_040} and \eqref{eq_062} with high probability, without loss of generality, we can assume that $\boldsymbol{X}$ has the property that $\boldsymbol{X}^\top\boldsymbol{U}_r$ is a positive semidefinite matrix, for this case, $\boldsymbol{U} = \boldsymbol{U}_r$.

\subsubsection{Proof of Theorem \ref{thm:main_psd}.}
\label{sec_proofs_main_step3}

In order to prove our main result, we need to first give a uniform upper bound of $K(\boldsymbol{X})$, then use the fact that for any local minimum $\widehat{\boldsymbol{X}}$, $K(\widehat{\boldsymbol{X}})\geqslant  0$, and finally solve for the range of possible $\widehat{\boldsymbol{X}}$. To start with, we first give an upper bound of perturbation terms. For simplicity of notations, denote $\nu_r \coloneqq   \|\boldsymbol{M}_r\|_{\ell_\infty}$. 

\begin{lemma}\label{lemma_k234}
If tuning parameters $\alpha,\lambda$ satisfy $100\sqrt{\nu_r} \leqslant \alpha\leqslant 200\sqrt{\nu_r}, 100 \|\boldsymbol{\Omega}-p\boldsymbol{J}\|\leqslant \lambda \leqslant 200 \|\boldsymbol{\Omega}-p\boldsymbol{J}\|$ and assume $ p\geqslant   C_1\frac{\log n}{n}$ with some absolute constant $C_1$. Then in an event $E$ with probability $\mathbb{P}[E]\geqslant  1-2n^{-3}$, uniformly for all $\boldsymbol{X}\in\mathbb{R}^{n\times r}$ and corresponding $\boldsymbol{\Delta}$ defined as before, we have
\begin{equation}\label{eq_115}
\begin{split}
	\sum_{i=2}^4 K_i(\boldsymbol{X}) \leqslant& 10^{-3}p\left[ \|\boldsymbol{\Delta}^\top\boldsymbol{\Delta}\|_F^2 + \|\boldsymbol{U}\boldsymbol{\Delta}^\top\|_F^2 \right]\\
 &  + C_2 p \sum_{i=1}^r \left\{\left[ C_3\left( \sqrt{\frac{n}{p}} +\frac{\log n}{p} \right)\nu_r + C_3 \sigma_{2r+1-i}- \sigma_i \right]_+\right\}^2 \\
 &+C_2p \frac{[(1-p)n+\log n/p]r\|\boldsymbol{N}\|_{\ell_\infty}^2}{p}   .
\end{split}
\end{equation}
\end{lemma}
Note in our proof of main theorem, we only use probabilistic tools in the above lemma to control perturbation terms, for the rest part of the proof, everything is deterministic.

Now denote $a \coloneqq \|\boldsymbol{\Delta}^\top\boldsymbol{\Delta}\|_F,b \coloneqq \|\boldsymbol{\Delta}^\top\boldsymbol{U}\|_F $ and 
\[
 \psi \coloneqq C_2 \left\{ \sum_{i=1}^r \left\{\left[ C_3\left( \sqrt{\frac{n}{p}} +\frac{\log n}{p} \right)\nu_r + C_3 \sigma_{2r+1-i}- \sigma_i \right]_+\right\}^2 + \frac{[(1-p)n+\log n/p]r\|\boldsymbol{N}\|_{\ell_\infty}^2}{p} \right\}.
\]

Putting Lemma \ref{lemma_k_psd} and Lemma \ref{lemma_k234} together and using the notations defined above we have
\begin{equation}\label{eq_080}
 \begin{split}
	 \frac{K(\boldsymbol{X})}{p} \leqslant & 1.001\|\boldsymbol{\Delta}\boldsymbol{\Delta}^\top\|_F^2 - 3\|\boldsymbol{X}\boldsymbol{X}^\top - \boldsymbol{U}\boldsymbol{U}^\top\|_F^2+10^{-3}\|\boldsymbol{U}\boldsymbol{\Delta}^\top\|_F^2 +\psi \\
	 =& 1.001a^2 - 3[\|\boldsymbol{\Delta}\boldsymbol{\Delta}^\top\|_F^2+2\|\boldsymbol{\Delta}\boldsymbol{U}^\top\|_F^2 + 2\langle \boldsymbol{\Delta}\boldsymbol{U}^\top,\boldsymbol{U} \boldsymbol{\Delta}^\top \rangle + 4\langle \boldsymbol{\Delta}\boldsymbol{\Delta}^\top, \boldsymbol{U}\boldsymbol{\Delta}^\top \rangle]\\
	 &  +10^{-3}\|\boldsymbol{U}\boldsymbol{\Delta}^\top\|_F^2+\psi,		
 \end{split}
\end{equation}
where second line use the decomposition 
\begin{equation}\label{eq_086}
 \begin{split}
	 \|\boldsymbol{X}\boldsymbol{X}^\top - \boldsymbol{U}\boldsymbol{U}^\top\|_F^2 =&\|\boldsymbol{U}\boldsymbol{\Delta}^\top + \boldsymbol{\Delta}\boldsymbol{U}^\top +\boldsymbol{\Delta}\boldsymbol{\Delta}^\top\|_F^2\\
	 = &  \|\boldsymbol{\Delta}\boldsymbol{\Delta}^\top\|_F^2+2\|\boldsymbol{\Delta}\boldsymbol{U}^\top\|_F^2 + 2\langle \boldsymbol{\Delta}\boldsymbol{U}^\top,\boldsymbol{U} \boldsymbol{\Delta}^\top \rangle + 4\langle \boldsymbol{\Delta}\boldsymbol{\Delta}^\top, \boldsymbol{U}\boldsymbol{\Delta}^\top \rangle.
 \end{split}
\end{equation}
By the definition of matrix inner product, we have
\begin{equation}\label{eq_082}
\begin{split}
 \|\boldsymbol{U}\boldsymbol{\Delta}^\top\|_F^2 =& \langle \boldsymbol{U}\boldsymbol{\Delta}^\top,\boldsymbol{U}\boldsymbol{\Delta}^\top\rangle= \operatorname{trace}(\boldsymbol{\Delta}\boldsymbol{U}^\top\boldsymbol{U}\boldsymbol{\Delta}^\top)= \operatorname{trace}(\boldsymbol{U}^\top\boldsymbol{U}\boldsymbol{\Delta}^\top\boldsymbol{\Delta})\\
	=& \langle \boldsymbol{U}^\top\boldsymbol{U},\boldsymbol{\Delta}^\top\boldsymbol{\Delta} \rangle,
\end{split}
\end{equation}
and
\begin{equation}\label{eq_083}
\langle \boldsymbol{\Delta}\boldsymbol{\Delta}^\top, \boldsymbol{U}\boldsymbol{\Delta}^\top\rangle =\operatorname{trace}(\boldsymbol{\Delta}\boldsymbol{\Delta}^\top\boldsymbol{U}\boldsymbol{\Delta}^\top) = \operatorname{trace}(\boldsymbol{\Delta}^\top\boldsymbol{\Delta}\boldsymbol{\Delta}^\top\boldsymbol{U})= \langle \boldsymbol{\Delta}^\top\boldsymbol{\Delta},\boldsymbol{\Delta}^\top\boldsymbol{U} \rangle.
\end{equation}
Moreover, since we choose $\boldsymbol{U}$ such that $\boldsymbol{U}^\top\boldsymbol{X}$ is positive semidefinite, $\boldsymbol{U}^\top\boldsymbol{\Delta}$ is a symmetric matrix and $\boldsymbol{U}^\top(\boldsymbol{\Delta}+\boldsymbol{U})$ is a positive semidefinite matrix. Therefore, we also have
\begin{equation}\label{eq_084}
\begin{split}
 \langle \boldsymbol{\Delta}\boldsymbol{U}^\top, \boldsymbol{U}\boldsymbol{\Delta}^\top \rangle =& \operatorname{trace}(\boldsymbol{U}\boldsymbol{\Delta}^\top\boldsymbol{U}\boldsymbol{\Delta}^\top)= \operatorname{trace}(\boldsymbol{\Delta}\boldsymbol{U}^\top\boldsymbol{\Delta}\boldsymbol{U}^\top) = \operatorname{trace}(\boldsymbol{U}^\top\boldsymbol{\Delta}\boldsymbol{U}^\top\boldsymbol{\Delta})\\
	=& \langle \boldsymbol{\Delta}^\top\boldsymbol{U},\boldsymbol{U}^\top\boldsymbol{\Delta} \rangle\\
	 =& \langle \boldsymbol{\Delta}^\top\boldsymbol{U},\boldsymbol{\Delta}^\top\boldsymbol{U}  \rangle\\
		=& \|\boldsymbol{\Delta}^\top\boldsymbol{U}\|_F^2
\end{split}
\end{equation}
and
\begin{equation}\label{eq_085}
 \langle \boldsymbol{\Delta}^\top\boldsymbol{\Delta}, \boldsymbol{U}^\top\boldsymbol{U} + \boldsymbol{\Delta}^\top \boldsymbol{U}\rangle = \langle \boldsymbol{\Delta}^\top\boldsymbol{\Delta},(\boldsymbol{U}+\boldsymbol{\Delta})^\top\boldsymbol{U}\rangle \geqslant  0.
\end{equation}
Here \eqref{eq_085} also uses the fact that inner product of two positive semidefinite matrices is non-negative. By putting \eqref{eq_080}, \eqref{eq_082}, \eqref{eq_083}, \eqref{eq_084}, \eqref{eq_085} together and using the notations defined before,
\begin{equation}\label{eq_101}
 \begin{split}
	 \frac{K(\boldsymbol{X})}{p} \leqslant & -1.999a^2 -  \langle \boldsymbol{\Delta}^\top\boldsymbol{\Delta}, 5.999\boldsymbol{U}^\top\boldsymbol{U} +12 \boldsymbol{\Delta}^\top \boldsymbol{U}\rangle- 6b^2  +\psi\\
	 \leqslant & -1.999a^2 - 6.001 \langle \boldsymbol{\Delta}^\top\boldsymbol{\Delta},\boldsymbol{\Delta}^\top \boldsymbol{U}  \rangle- 6b^2  +\psi\\
	 \leqslant &-1.999 a^2 + 6.001ab-6b^2 +\psi
 \end{split}
\end{equation}
holds for all $\boldsymbol{X}\in\mathbb{R}^{n\times r}$.
For the last line, we apply Cauchy-Schwarz inequality for matrices again.

Note for any local minimum $\widehat{\boldsymbol{X}}$, we have $K(\widehat{\boldsymbol{X}})\geqslant  0$. Combining with \eqref{eq_101} we have
\[
 -1.999 a^2 + 6.001ab-6b^2 +\psi \geqslant  0.
\]
Solving for $a$ and $b$ we have 
\begin{equation}\label{eq_087}
 a\leqslant C_4 \sqrt{\psi},\qquad  b\leqslant  C_4\sqrt{\psi},
\end{equation}
here we also use the fact that $a,b\geqslant  0$.

Now use the fact that $K(\widehat{\boldsymbol{X}})\geqslant  0$ again together with \eqref{eq_080}, \eqref{eq_086} and \eqref{eq_087},
\begin{equation}\label{eq_116}
\begin{split}
 & 2.999\|\widehat{\boldsymbol{X}}\widehat{\boldsymbol{X}}^\top - \boldsymbol{U}\boldsymbol{U}^\top\|_F^2\\
	\leqslant & 1.001\|\boldsymbol{\Delta}\boldsymbol{\Delta}^\top\|_F^2 - 10^{-3}\|\widehat{\boldsymbol{X}}\widehat{\boldsymbol{X}}^\top - \boldsymbol{U}\boldsymbol{U}^\top\|_F^2 + 10^{-3}\|\boldsymbol{U}\boldsymbol{\Delta}^\top\|_F^2+\psi\\
 \leqslant & a^2 +\psi +4\times 10^{-3}ab\\
	\leqslant &C_4\psi.
\end{split}
\end{equation}
Therefore,
\begin{equation*}	
 \begin{split}
	 \|\widehat{\boldsymbol{X}}\widehat{\boldsymbol{X}}^\top - \boldsymbol{U}\boldsymbol{U}^\top\|_F^2 \leqslant &    C_2  \sum_{i=1}^r \left\{\left[ C_3\left( \sqrt{\frac{n}{p}} +\frac{\log n}{p} \right)\nu_r + C_3 \sigma_{2r+1-i}- \sigma_i \right]_+\right\}^2\\
	 & +C_2 \frac{[(1-p)n+\log n/p]r\|\boldsymbol{N}\|_{\ell_\infty}^2}{p} 
 \end{split}
\end{equation*}
and also 
\begin{equation*}
 \begin{split}
	 \|\widehat{\boldsymbol{X}}\widehat{\boldsymbol{X}}^\top-\boldsymbol{M} \|_F^2 \leqslant& C_2  \sum_{i=1}^r \left\{\left[ C_3\left( \sqrt{\frac{n}{p}} +\frac{\log n}{p} \right)\nu_r+ C_3 \sigma_{2r+1-i} - \sigma_i \right]_+\right\}^2\\
	 & +C_2 \frac{[(1-p)n+\log n/p]r\|\boldsymbol{N}\|_{\ell_\infty}^2}{p}  +\|\boldsymbol{N}\|_F^2. 
 \end{split}
\end{equation*}
Here we use the fact that 
$$\|\widehat{\boldsymbol{X}}\widehat{\boldsymbol{X}}^\top-\boldsymbol{M} \|_F^2 = \|\widehat{\boldsymbol{X}}\widehat{\boldsymbol{X}}^\top - \boldsymbol{U}\boldsymbol{U}^\top\|_F^2 -2 \langle\widehat{\boldsymbol{X}}\widehat{\boldsymbol{X}}^\top,\boldsymbol{N}  \rangle + \|\boldsymbol{N}\|_F^2 \leqslant \|\widehat{\boldsymbol{X}}\widehat{\boldsymbol{X}}^\top - \boldsymbol{U}\boldsymbol{U}^\top\|_F^2 + \|\boldsymbol{N}\|_F^2$$
and the inequality holds since the inner product of two positive semidefinite matrices is non-negative.


\subsection{Proof of lemmas}
\label{sec_proofs_proof}

Here we present proof for lemmas used in former sections.

\subsubsection{A proof of Lemma \ref{lemma_concern}}
\begin{proof}

 First of all, by using the definition of matrix inner product and Hadamard product, we have
 \begin{equation}\label{eq_090}
 \begin{split}
	 | \langle \mathcal{P}_{\Omega_0}(\boldsymbol{A}\boldsymbol{C}^\top),\mathcal{P}_{\Omega_0}(\boldsymbol{B}\boldsymbol{D}^\top)\rangle-t  \langle \boldsymbol{A}\boldsymbol{C}^\top,\boldsymbol{B}\boldsymbol{D}^\top \rangle|  =& |\langle \boldsymbol{\Omega}_0-t\boldsymbol{J}, (\boldsymbol{A}\boldsymbol{C}^\top \circ \boldsymbol{B}\boldsymbol{D}^\top)\rangle|\\
		\leqslant& \|\boldsymbol{\Omega}_0-t\boldsymbol{J}\| \|(\boldsymbol{A}\boldsymbol{C}^\top \circ \boldsymbol{B}\boldsymbol{D}^\top) \|_*, 
	 \end{split}
 \end{equation}
 where the inequality holds by matrix H\"{o}lder's inequality.
 So the only thing left over is to give a bound of $\|(\boldsymbol{A}\boldsymbol{C}^\top \circ \boldsymbol{B}\boldsymbol{D}^\top) \|_*$. Notice one can decompose the matrix into sum of rank one matrices as following
 \begin{equation*}
	 \begin{split}
			\boldsymbol{A}\boldsymbol{C}^\top \circ \boldsymbol{B}\boldsymbol{D}^\top =& \left(\sum_{k=1}^{r_1} \boldsymbol{A}_{\cdot,k}\boldsymbol{C}_{\cdot,k}^\top\right)\circ \left(\sum_{k=1}^{r_2} \boldsymbol{B}_{\cdot,k}\boldsymbol{D}_{\cdot,k}^\top\right) = \sum_{l=1}^{r_1}\sum_{m=1}^{r_2} (\boldsymbol{A}_{\cdot,l}\circ \boldsymbol{B}_{\cdot,m}) (\boldsymbol{C}_{\cdot,l}\circ \boldsymbol{D}_{\cdot,m})^\top.
	 \end{split}
 \end{equation*}
 So one can upper bound the nuclear norm via
 \begin{equation*}
	 \begin{split}
	 \|(\boldsymbol{A}\boldsymbol{C}^\top \circ \boldsymbol{B}\boldsymbol{D}^\top) \|_* \leqslant &  \sum_{l=1}^{r_1}\sum_{m=1}^{r_2}\|(\boldsymbol{A}_{\cdot,l}\circ \boldsymbol{B}_{\cdot,m}) (\boldsymbol{C}_{\cdot,l}\circ \boldsymbol{D}_{\cdot,m})^\top \|_*\\
	 =& \sum_{l=1}^{r_1}\sum_{m=1}^{r_2}\| \boldsymbol{A}_{\cdot,l}\circ \boldsymbol{B}_{\cdot,m}\|_2\|\boldsymbol{C}_{\cdot,l}\circ \boldsymbol{D}_{\cdot,m} \|_2\\
	 =& \sum_{l=1}^{r_1}\sum_{m=1}^{r_2} \sqrt{\sum_{k=1}^{n_1} A_{k,l}^2 B_{k,m}^2}\sqrt{\sum_{k=1}^{n_2} C_{k,l}^2 D_{k,m}^2},
	 \end{split}
 \end{equation*}
 where first line is by triangle inequality and we can replace nuclear norm by vector $\ell_2$ norms in second line since summands are all rank one matrices. Now apply Cauchy-Schwarz inequality twice we have
 \begin{equation}\label{eq_091}
	 \begin{split}
	 \|(\boldsymbol{A}\boldsymbol{C}^\top \circ \boldsymbol{B}\boldsymbol{D}^\top) \|_* 	\leqslant & \sqrt{\sum_{l=1}^{r_1}\sum_{m=1}^{r_2}\sum_{k=1}^{n_1} A_{k,l}^2 B_{k,m}^2}\sqrt{\sum_{l=1}^{r_1}\sum_{m=1}^{r_2}\sum_{k=1}^{n_2} C_{k,l}^2 D_{k,m}^2}\\
	 =& \sqrt{\sum_{k=1}^{n_1}\|\boldsymbol{A}_{k,\cdot}\|_2^2 \|\boldsymbol{B}_{k,\cdot}\|_2^2}\sqrt{\sum_{k=1}^{n_2}\|\boldsymbol{C}_{k,\cdot}\|_2^2 \|\boldsymbol{D}_{k,\cdot}\|_2^2}.
	 \end{split}
 \end{equation}

 Putting \eqref{eq_090} and \eqref{eq_091} together we have
 \begin{equation*}
	 \begin{split}
		&| \langle \mathcal{P}_{\Omega_0}(\boldsymbol{A}\boldsymbol{C}^\top),\mathcal{P}_{\Omega_0}(\boldsymbol{B}\boldsymbol{D}^\top)\rangle-t  \langle \boldsymbol{A}\boldsymbol{C}^\top,\boldsymbol{B}\boldsymbol{D}^\top \rangle|\\
			\leqslant& \|\boldsymbol{\Omega}_0-t\boldsymbol{J}\| \sqrt{\sum_{k=1}^{n_1}\|\boldsymbol{A}_{k,\cdot}\|_2^2 \|\boldsymbol{B}_{k,\cdot}\|_2^2}\sqrt{\sum_{k=1}^{n_2}\|\boldsymbol{C}_{k,\cdot}\|_2^2 \|\boldsymbol{D}_{k,\cdot}\|_2^2} .
	 \end{split}
 \end{equation*}
\end{proof}

\subsubsection{A proof of Lemma \ref{lemma_k234}}
\begin{proof}

The proof of Lemma \ref{lemma_k234} can be divided into controlling $K_2(\boldsymbol{X})$, $K_3(\boldsymbol{X})$ and $K_4(\boldsymbol{X})$.

For $K_2(\boldsymbol{X})$, we have

\begin{lemma}\label{lemma_k2_psd}
In an event $E_a$ with probability $\mathbb{P}[E_a]\geqslant  1-n^{-3}$, uniformly for all $\boldsymbol{X}\in\mathbb{R}^{n\times r}$ and corresponding $\boldsymbol{\Delta}$ defined as before, we have
 \begin{equation*}
	 K_2(\boldsymbol{X}) \leqslant  \|\boldsymbol{\Omega}-p\boldsymbol{J}\| \left[ 19\sum_{i=1}^n \|\boldsymbol{\Delta}_{i,\cdot}\|_2^4 + 18\nu_r  \|\boldsymbol{\Delta}\|_F^2 + 9\nu_r  \sum_{i=s+1}^r \sigma_i \right] +3\times 10^{-4}p \|\boldsymbol{U}\boldsymbol{\Delta}^\top\|_F^2,
		 \end{equation*}
		 where $s$ is defined by
		 \begin{equation}\label{eq_104}
			 s \coloneqq \max \left\{s\leqslant r ,\;\sigma_s\geqslant  C_p \frac{\nu_r \log n}{p}\right\}
		 \end{equation}
		 with $C_p$ an absolute constant. Set $s = 0$ if $\sigma_1<C_p \frac{\nu_r \log n}{p}$.
\end{lemma}

For $K_3(\boldsymbol{X})$, we use a modified version of \citet[Lemma 11]{ge2017no}:
\begin{lemma}[{\citealt[Lemma 11]{ge2017no}}]\label{lemma_k3_psd}
 If $\alpha \geqslant  100\sqrt{\nu_r}$, then uniformly for all $\boldsymbol{X}\in\mathbb{R}^{n\times r}$ and corresponding $\boldsymbol{\Delta}$ defined as before, we have
 \begin{equation*}
	 K_3(\boldsymbol{X}) \leqslant 199.54\lambda\alpha^2 \|\boldsymbol{\Delta}\|_F^2 - 0.3\lambda \sum_{i=1}^n \|\boldsymbol{\Delta}_{i,\cdot}\|_2^4.
 \end{equation*}
\end{lemma}
The main difference is that we keep the extra negative term. We will give a proof in appendix for completeness.

For $K_4(\boldsymbol{X})$, we have

\begin{lemma}\label{lemma_k4_psd}
 Uniformly for all $\boldsymbol{X}\in\mathbb{R}^{n\times r}$ and corresponding $\boldsymbol{\Delta}$ defined as before, we have 
	 \begin{equation*}
 \begin{split}
	 K_4(\boldsymbol{X})\leqslant& 5\times 10^{-4}p \|\boldsymbol{\Delta}\boldsymbol{\Delta}^\top\|_F^2 +2\times 10^{-4} p \|\boldsymbol{U}\boldsymbol{\Delta}^\top\|_F^2  + C_2 \frac{r\| \mathcal{P}_{\Omega} (\boldsymbol{N}) - p\boldsymbol{N} \|^2}{p}\\
	 & + 6p\langle \boldsymbol{\Delta}\boldsymbol{\Delta}^\top ,\boldsymbol{N} \rangle.
	 \end{split}
 \end{equation*}
\end{lemma}

Now we can use Lemma \ref{lemma_vu14_psd} together with Lemma \ref{lemma_spec} to bound $\| \mathcal{P}_{\Omega} (\boldsymbol{N}) - p\boldsymbol{N} \|$ and $\|\boldsymbol{\Omega}-p\boldsymbol{J}\|$ (similar result can also be found in \citet{keshavan2010matrix}): As long as $  p\geqslant C_1 \frac{\log n}{n}$ with some absolute constant $C_1$, they are bounded by 
$$\left(C\sqrt{np(1-p)}+C\sqrt{\log n}\right)\|\boldsymbol{N}\|_{\ell_\infty}$$ and $C\sqrt{np}$ correspondingly in an event $E_b$ with probability $\mathbb{P}[E_b]\geqslant  1-n^{-3}$ .

Note that we choose $\alpha,\lambda$ such that $100\sqrt{\nu_r}\leqslant \alpha\leqslant 200\sqrt{\nu_r}, 100\|\boldsymbol{\Omega}-p\boldsymbol{J}\|\leqslant \lambda \leqslant 200\|\boldsymbol{\Omega}-p\boldsymbol{J}\|$. By using union bound to put the estimates together,
\begin{equation}\label{eq_102}
 \begin{split}
\sum_{i=2}^4 K_i(\boldsymbol{X}) \leqslant& 5\times 10^{-4}p\left[ \|\boldsymbol{\Delta}^\top\boldsymbol{\Delta}\|_F^2 + \|\boldsymbol{U}\boldsymbol{\Delta}^\top\|_F^2 \right] +C_2 [(1-p)n+\log n/p]r\|\boldsymbol{N}\|_{\ell_\infty}^2\\
&   + C_3\sqrt{np}\nu_r \|\boldsymbol{\Delta}\|_F^2  + C_4 \sqrt{np} \nu_r \sum_{i=s+1}^r \sigma_i + 6p\langle \boldsymbol{\Delta}\boldsymbol{\Delta}^\top ,\boldsymbol{N} \rangle,
 \end{split}
\end{equation}
holds in an event $E$ with probability $\mathbb{P}[E]\geqslant  1-2n^{-3}$. 

For $\|\boldsymbol{\Delta}^\top\boldsymbol{\Delta}\|_F^2$, we have  
\begin{equation}\label{eq_093}
 \|\boldsymbol{\Delta}^\top\boldsymbol{\Delta}\|_F^2  = \langle \boldsymbol{\Delta}^\top\boldsymbol{\Delta},\boldsymbol{\Delta}^\top\boldsymbol{\Delta} \rangle = \sum_{i=1}^r \sigma_i^4(\boldsymbol{\Delta}),
\end{equation}
where $\sigma_i(\boldsymbol{\Delta})$ denotes $i$-th largest singular value of $\boldsymbol{\Delta}$.

In order to proceed, we need the following result:
\begin{lemma}[{\citealt[Problem III.6.14]{bhatia2013matrix}}] \label{claim:01}
 Let $\boldsymbol{A},\boldsymbol{B}\in \mathbb{R}^{n\times n}$ be two symmetric matrices, $\lambda_1(\boldsymbol{A})\geqslant  \lambda_2(\boldsymbol{A})\geqslant  \cdots \geqslant  \lambda_n(\boldsymbol{A})$ and $\lambda_1(\boldsymbol{B})\geqslant  \lambda_2(\boldsymbol{B})\geqslant  \cdots \geqslant  \lambda_n(\boldsymbol{B})$ are eigenvalues of $\boldsymbol{A}$ and $\boldsymbol{B}$. Then the following holds:
	 \[
		 \sum_{i=1}^n \lambda_{i}(\boldsymbol{A})\lambda_{n+1-i}(\boldsymbol{B})\leqslant  \langle \boldsymbol{A},\boldsymbol{B} \rangle \leqslant \sum_{i = 1}^n \lambda_i(\boldsymbol{A}) \lambda_i(\boldsymbol{B}).
	 \] 
\end{lemma} 
This result can also be derived from Schur-Horn theorem (see, e.g., \citet[Theorem 9.B.1, Theorem 9.B.2]{marshall2011inequalities}) together with Abel's summation formula.

From Lemma \ref{claim:01}, we have
\begin{equation}\label{eq_094}
 \begin{split}
	 \|\boldsymbol{U}\boldsymbol{\Delta}^\top\|_F^2 =&  \operatorname{trace}(\boldsymbol{\Delta}\boldsymbol{U}^\top\boldsymbol{U}\boldsymbol{\Delta}^\top)\\
		= & \langle \boldsymbol{U}^\top\boldsymbol{U}, \boldsymbol{\Delta}^\top \boldsymbol{\Delta}  \rangle \\
		\geqslant  & \sum_{i=1}^r \lambda_{r+1-i}(\boldsymbol{U}^\top \boldsymbol{U})\lambda_i(\boldsymbol{\Delta}^\top\boldsymbol{\Delta})\\
		= & \sum_{i=1}^r \sigma_i^2(\boldsymbol{\Delta})\sigma_{r+1-i}^2(\boldsymbol{U}),
 \end{split}
 \end{equation}
 and 
 \begin{equation}\label{eq_117}
	 \begin{split}
		 \langle \boldsymbol{\Delta}\boldsymbol{\Delta}^\top, \boldsymbol{N} \rangle \leqslant&  \sum_{i=1}^n \lambda_i(\boldsymbol{\Delta}\boldsymbol{\Delta}^\top) \lambda_i(\boldsymbol{N})\\
		 = & \sum_{i=1}^r \sigma_i^2(\boldsymbol{\Delta})\sigma_i(\boldsymbol{N}).
	 \end{split}
 \end{equation}
 Here we use the fact that $\lambda_i(\boldsymbol{U}^\top\boldsymbol{U}) = \sigma_i^2(\boldsymbol{U})$, $\lambda_i(\boldsymbol{\Delta}^\top\boldsymbol{\Delta}) = \sigma_i^2(\boldsymbol{\Delta})$, $\lambda_i(\boldsymbol{N}) = \sigma_i(\boldsymbol{N})$ and 
 \[
	 \lambda_i(\boldsymbol{\Delta}\boldsymbol{\Delta}^\top) = \left\{ \begin{array}{ll}
			\sigma_i^2(\boldsymbol{\Delta}) & i = 1,\cdots,r\\
			0 & i = r+1,\cdots,n.
	 \end{array} \right.
 \]

Putting \eqref{eq_093}, \eqref{eq_094} and \eqref{eq_117} together we have 
\begin{equation}\label{eq_103}
 \begin{split}
& - 5\times 10^{-4}p\left[ \|\boldsymbol{\Delta}^\top\boldsymbol{\Delta}\|_F^2 + \|\boldsymbol{U}\boldsymbol{\Delta}^\top\|_F^2 \right] + C_3\sqrt{np}\nu_r \|\boldsymbol{\Delta}\|_F^2 + 6p\langle \boldsymbol{\Delta}\boldsymbol{\Delta}^\top, \boldsymbol{N} \rangle\\
 \leqslant  &   5\times 10^{-4} p\sum_{i=1}^r \left\{ - \sigma_i^4(\boldsymbol{\Delta})+\left[ C_3 \sqrt{\frac{n }{p}}\nu_r - \sigma_{r+1-i}^2(\boldsymbol{U})+C_3 \sigma_i(\boldsymbol{N}) \right] \sigma_i^2(\boldsymbol{\Delta})  \right\}\\
 \leqslant & C_2 p \sum_{i=1}^r\left\{\left[ C_3 \sqrt{\frac{n }{p}}\nu_r + C_3 \sigma_{2r+1-i}- \sigma_i \right]_+\right\}^2 .
 \end{split}
\end{equation}
Here we optimize over a series of quadratic functions of $\sigma_i^2(\boldsymbol{\Delta})$ in the last line. In the last line, we also use the fact that $\sigma_i(\boldsymbol{N}) = \sigma_{r+i}(\boldsymbol{M}) = \sigma_{r+i}$. Finally putting \eqref{eq_102} and \eqref{eq_103} together we have

\begin{equation*}
\begin{split} 
	\sum_{i=2}^4 K_i(\boldsymbol{X}) \leqslant & 10^{-3}p\left[ \|\boldsymbol{\Delta}^\top\boldsymbol{\Delta}\|_F^2 + \|\boldsymbol{U}\boldsymbol{\Delta}^\top\|_F^2 \right] +C_2 \left[(1-p)n+\log n/p\right]r\|\boldsymbol{N}\|_{\ell_\infty}^2\\
 &C_2 p \sum_{i=1}^r\left\{\left[ C_3 \sqrt{\frac{n }{p}}\nu_r + C_3 \sigma_{2r+1-i}- \sigma_i \right]_+\right\}^2  + C_4 \sqrt{np } \nu_r \sum_{i=s+1}^r \sigma_i\\
 \leqslant &  10^{-3}p\left[ \|\boldsymbol{\Delta}^\top\boldsymbol{\Delta}\|_F^2 + \|\boldsymbol{U}\boldsymbol{\Delta}^\top\|_F^2 \right]\\
 & + C_2 p \sum_{i=1}^r \left\{\left[ C_3\left( \sqrt{\frac{n }{p}} +\frac{\log n}{p} \right)\nu_r + C_3 \sigma_{2r+1-i}- \sigma_i \right]_+\right\}^2\\
 & + C_2p \frac{\left[(1-p)n + \log n/p\right]r\|\boldsymbol{N}\|_{\ell_\infty}^2}{p}  ,
\end{split}
\end{equation*}
where the last inequality holds since by definition of $s$, for any $i>s$, we have $\sigma_i < C_p \frac{\nu_r\log n}{p}$. Choosing $C_3$ sufficient large we have
\begin{equation*}
 \begin{split}
 \left\{\left[ C_3\left( \sqrt{\frac{n }{p}}+\frac{\log n}{p}   \right)\nu_r + C_3 \sigma_{2r+1-i}-\sigma_i \right]_+\right\}^2 \geqslant&  \left[C_3\sqrt{\frac{n }{p}}\nu_r +\sigma_i \right]^2 \geqslant    \sqrt{\frac{n }{p}}\nu_r \sigma_i ,		
 \end{split}
\end{equation*}
which finishes the proof. 
\end{proof}

\subsubsection{A proof of Lemma \ref{lemma_k2_psd}}

\begin{proof}

Recall that we define $\boldsymbol{\Delta}$ as $\boldsymbol{\Delta}\coloneqq \boldsymbol{X}-\boldsymbol{U}$, $D_{\Omega,p}(\boldsymbol{X}\boldsymbol{X}^\top-\boldsymbol{U}\boldsymbol{U}^\top,\boldsymbol{X}\boldsymbol{X}^\top-\boldsymbol{U}\boldsymbol{U}^\top)$ can be decomposed as following
\begin{equation}\label{eq_031}
 \begin{split}
		&D_{\Omega,p}(\boldsymbol{X}\boldsymbol{X}^\top-\boldsymbol{U}\boldsymbol{U}^\top,\boldsymbol{X}\boldsymbol{X}^\top-\boldsymbol{U}\boldsymbol{U}^\top)\\
		=& D_{\Omega,p}(\boldsymbol{U}\boldsymbol{\Delta}^\top +\boldsymbol{\Delta}\boldsymbol{U}^\top+\boldsymbol{\Delta}\boldsymbol{\Delta}^\top,\boldsymbol{U}\boldsymbol{\Delta}^\top +\boldsymbol{\Delta}\boldsymbol{U}^\top+\boldsymbol{\Delta}\boldsymbol{\Delta}^\top) \\
	 =& D_{\Omega,p}(\boldsymbol{U}\boldsymbol{\Delta}^\top+\boldsymbol{\Delta}\boldsymbol{U}^\top,\boldsymbol{U}\boldsymbol{\Delta}^\top+\boldsymbol{\Delta}\boldsymbol{U}^\top)+D_{\Omega,p}(\boldsymbol{\Delta}\boldsymbol{\Delta}^\top,\boldsymbol{\Delta}\boldsymbol{\Delta}^\top)  +  4D_{\Omega,p}(\boldsymbol{U}\boldsymbol{\Delta}^\top,\boldsymbol{\Delta}\boldsymbol{\Delta}^\top).
 \end{split}
\end{equation}
Here we use the fact that $\Omega$ is symmetric. Our strategy here is using Lemma \ref{lemma_rip_psd} to give a tight bound to as many as possible terms, for those terms that Lemma \ref{lemma_rip_psd} cannot handle, we use Lemma \ref{lemma_concern} to give a bound. First for the second and third term of \eqref{eq_031}, as Lemma \ref{lemma_rip_psd} cannot apply here, we use Lemma \ref{lemma_concern} to give

\begin{equation}\label{eq_032}
|D_{\Omega,p}(\boldsymbol{\Delta}\boldsymbol{\Delta}^\top,\boldsymbol{\Delta}\boldsymbol{\Delta}^\top)|\leqslant  \|\boldsymbol{\Omega}-p\boldsymbol{J}\|\sum_{i=1}^n \|\boldsymbol{\Delta}_{i,\cdot}\|_2^4
\end{equation}
and
\begin{equation}\label{eq_061}
\begin{split}
 4|D_{\Omega,p}(\boldsymbol{U}\boldsymbol{\Delta}^\top,\boldsymbol{\Delta}\boldsymbol{\Delta}^\top)|\leqslant & 4\|\boldsymbol{\Omega}-p\boldsymbol{J}\|\sqrt{\sum_{i=1}^n \|\boldsymbol{U}_{i,\cdot}\|_2^2 \|\boldsymbol{\Delta}_{i,\cdot}\|_2^2}\sqrt{\sum_{i=1}^n \|\boldsymbol{\Delta}_{i,\cdot}\|_2^4}\\
	\leqslant &2\|\boldsymbol{\Omega}-p\boldsymbol{J}\|\nu_r \|\boldsymbol{\Delta}\|_F^2  +2\|\boldsymbol{\Omega}-p\boldsymbol{J}\|\sum_{i=1}^n \|\boldsymbol{\Delta}_{i,\cdot}\|_2^4, 
\end{split}
\end{equation}
where for the second inequality we use the fact that $2xy\leqslant x^2+y^2$.

Finally for the first term of \eqref{eq_031}, if $\boldsymbol{U}$ is good enough such that the incoherence $\mu(\boldsymbol{U})$ is well-bounded, then we can apply Lemma \ref{lemma_rip_psd} directly and get a tight bound. If $\mu(\boldsymbol{U})$ is not good enough, we want to split $\boldsymbol{U}$ into two parts and hope first few columns have good incoherence. To be more precise, recall that we assume $\boldsymbol{U} = \boldsymbol{U}_r = [\sqrt{\sigma_1}\boldsymbol{u}_1\; \dots\; \sqrt{\sigma_r}\boldsymbol{u}_r]$, similar to \eqref{eq_100}, for the incoherence of the first $k$ columns, we have 
\[
	 \begin{split}
		 \mu\left( \colspan([\sqrt{\sigma_1}\boldsymbol{u}_1\; \dots\; \sqrt{\sigma_k}\boldsymbol{u}_k]) \right)  = & \frac{n}{k}\max_i \sum_{j = 1}^{k} u_{i,j}^2\\
		 \leqslant & \frac{n}{k\sigma_k}\max_i \sum_{j = 1}^{k} \sigma_j u_{i,j}^2\\
		 \leqslant & \frac{n}{k\sigma_k}\max_i \sum_{j = 1}^{r} \sigma_j u_{i,j}^2\\
		 \leqslant& \frac{n\nu_r}{k\sigma_k},
	 \end{split}
\]
where $\mu(\cdot)$ is defined in \eqref{eq:mu_subspace}.

For fixed $s$ defined as in \eqref{eq_104}, denote first $s$ columns of $\boldsymbol{U}$ as $\boldsymbol{U}^{1}$, and remaining part as $\boldsymbol{U}^{2}$. Then $\mu\left(\colspan(\boldsymbol{U}^{1})\right) \leqslant \frac{n \nu_r  }{s\sigma_s}$, which makes it possible to apply Lemma \ref{lemma_rip_psd} to space spanned by $\boldsymbol{U}^{1}$ since we have $\sigma_s\geqslant  C_p \frac{\nu_r \log n}{p}$. Decompose $\boldsymbol{U}$ as $\boldsymbol{U} = [\boldsymbol{U}^{1}\; \boldsymbol{U}^{2}]$, and $\boldsymbol{\Delta} $ can also be decomposed as $\boldsymbol{\Delta}  = [\boldsymbol{\Delta}^{1}\; \boldsymbol{\Delta}^{2}]$ correspondingly. Note by our assumption that $\boldsymbol{U} = \boldsymbol{U}_r$, we have $(\boldsymbol{U}^{1})^\top \boldsymbol{U}^{2} = \boldsymbol{0}$. So we can further decompose the first term of \eqref{eq_031} as

\begin{equation}\label{eq_030}
\begin{split}
	&D_{\Omega,p}(\boldsymbol{U}\boldsymbol{\Delta}^\top+\boldsymbol{\Delta}\boldsymbol{U}^\top,\boldsymbol{U}\boldsymbol{\Delta}^\top+\boldsymbol{\Delta}\boldsymbol{U}^\top)\\
	 =&  D_{\Omega,p}\left([\boldsymbol{U}^{1}\; \boldsymbol{U}^{2}]  [\boldsymbol{\Delta}^{1}\; \boldsymbol{\Delta}^{2}]^\top+  [\boldsymbol{\Delta}^{1}\; \boldsymbol{\Delta}^{2}][\boldsymbol{U}^{1}\; \boldsymbol{U}^{2}]^\top,[\boldsymbol{U}^{1}\; \boldsymbol{U}^{2}]  [\boldsymbol{\Delta}^{1}\; \boldsymbol{\Delta}^{2}]^\top\right.\\
	 &\left.\qquad+  [\boldsymbol{\Delta}^{1}\; \boldsymbol{\Delta}^{2}][\boldsymbol{U}^{1}\; \boldsymbol{U}^{2}]^\top \right) \\
 = & D_{\Omega,p}\left(\boldsymbol{U}^{1}(\boldsymbol{\Delta}^{1})^\top + \boldsymbol{\Delta}^{1}(\boldsymbol{U}^{1})^\top,\boldsymbol{U}^{1}(\boldsymbol{\Delta}^{1})^\top + \boldsymbol{\Delta}^{1}(\boldsymbol{U}^{1})^\top\right)   \\
 & +  4 D_{\Omega,p}\left(\boldsymbol{U}^{1}(\boldsymbol{\Delta}^{1})^\top,\boldsymbol{U}^{2}(\boldsymbol{\Delta}^{2})^\top\right) + 2D_{\Omega,p}\left(\boldsymbol{U}^{2}(\boldsymbol{\Delta}^{2})^\top,\boldsymbol{U}^{2}(\boldsymbol{\Delta}^{2})^\top\right)   \\
 &+  2D_{\Omega,p}\left(\boldsymbol{U}^{2}(\boldsymbol{\Delta}^{2})^\top,\boldsymbol{\Delta}^{2}(\boldsymbol{U}^{2})^\top\right)+ 4D_{\Omega,p}\left(\boldsymbol{U}^{1}(\boldsymbol{\Delta}^{1})^\top,\boldsymbol{\Delta}^{2}(\boldsymbol{U}^{2})^\top\right) .
\end{split}
\end{equation}

Now we can apply tight approximation Lemma \ref{lemma_rip_psd} to the first term of \eqref{eq_030}. If we choose $C_p$ sufficient large such that for $s$ defined as before, $p \geqslant  C \frac{\nu_r  \log n}{\delta^2 \sigma_s}\geqslant  C\frac{\mu\left(\colspan(\boldsymbol{U}^{1})\right) s \log n}{\delta^2 n}$ with $\delta = 2.5\times 10^{-5}$, then 
\begin{equation}\label{eq_105}
\begin{split}
&\left| D_{\Omega,p}\left(\boldsymbol{U}^{1}(\boldsymbol{\Delta}^{1})^\top + \boldsymbol{\Delta}^{1}(\boldsymbol{U}^{1})^\top,\boldsymbol{U}^{1}(\boldsymbol{\Delta}^{1})^\top + \boldsymbol{\Delta}^{1}(\boldsymbol{U}^{1})^\top\right)\right| \\
 \leqslant& 2.5\times 10^{-5}p \|\boldsymbol{U}^{1}(\boldsymbol{\Delta}^{1})^\top+ \boldsymbol{\Delta}^{1}(\boldsymbol{U}^{1})^\top\|_F^2\\
 \leqslant & 5\times 10^{-5}p (\|\boldsymbol{U}^{1}(\boldsymbol{\Delta}^{1})^\top\|_F^2 +\|\boldsymbol{\Delta}^{1}(\boldsymbol{U}^{1})^\top\| _F^2) \\
 \leqslant & 10^{-4}p\|\boldsymbol{U}\boldsymbol{\Delta}^\top\|_F^2
\end{split}
\end{equation}
holds in an event $E_a$ with probability $\mathbb{P}[E_a]\geqslant  1-n^{-3}$, where the second inequality uses the fact that $(x+y)^2\leqslant 2x^2+2y^2$, and last inequality uses the fact that $(\boldsymbol{U}^{1})^\top \boldsymbol{U}^{2} = \boldsymbol{0}$.

For the rest terms in \eqref{eq_030}, by applying Lemma \ref{lemma_concern} we have
\begin{equation}\label{eq_106}
\begin{split}
 4|  D_{\Omega,p}(\boldsymbol{U}^{1}(\boldsymbol{\Delta}^{1})^\top,\boldsymbol{U}^{2}(\boldsymbol{\Delta}^{2})^\top)| \leqslant &4 \|\boldsymbol{\Omega}-p\boldsymbol{J}\| \sqrt{\sum_{i=1}^n \|\boldsymbol{U}_{i,\cdot}^{1}\|_2^2\|\boldsymbol{U}_{i,\cdot}^{2}\|_2^2} \sqrt{\sum_{i=1}^n \|\boldsymbol{\Delta}_{i,\cdot}^{1}\|_2^2\|\boldsymbol{\Delta}_{i,\cdot}^{2}\|_2^2}\\
	\leqslant & 2\|\boldsymbol{\Omega}-p\boldsymbol{J}\|\left[ \nu_r  \|\boldsymbol{U}^{2}\|_F^2 +\sum_{i=1}^n \|\boldsymbol{\Delta}_{i,\cdot}\|_2^4 \right]
\end{split}
\end{equation}
for the second term in \eqref{eq_030}, where the second inequality use the fact that $\|\boldsymbol{U}_{i,\cdot}^{1}\|_2^2 \leqslant \|\boldsymbol{U}_{i,\cdot} \|_2^2\leqslant \nu_r,\|\boldsymbol{\Delta}_{i,\cdot}^{1}\|_2^2 \leqslant \|\boldsymbol{\Delta}_{i,\cdot}\|_2^2, \|\boldsymbol{\Delta}_{i,\cdot}^{2}\|_2^2 \leqslant \|\boldsymbol{\Delta}_{i,\cdot}\|_2^2$ and $2xy\leqslant x^2+y^2$. For the third term, applying Lemma \ref{lemma_concern} again we have
\begin{equation}\label{eq_107}
\begin{split}
2 | D_{\Omega,p}(\boldsymbol{U}^{2}(\boldsymbol{\Delta}^{2})^\top,\boldsymbol{U}^{2}(\boldsymbol{\Delta}^{2})^\top) | \leqslant & 2\|\boldsymbol{\Omega}-p\boldsymbol{J}\|\sqrt{\sum_{i=1}^n \|\boldsymbol{U}_{i,\cdot}^{2}\|_2^4}\sqrt{\sum_{i=1}^n \|\boldsymbol{\Delta}_{i,\cdot}^{2}\|_2^4}\\
 \leqslant &  \|\boldsymbol{\Omega}-p\boldsymbol{J}\|	\left[ \nu_r  \|\boldsymbol{U}^{2}\|_F^2 +\sum_{i=1}^n \|\boldsymbol{\Delta}_{i,\cdot}\|_2^4 \right],
\end{split}
\end{equation}
where for the second inequality we also use the properties used in bounding second term. For the fourth and last term in \eqref{eq_030}, applying Lemma \ref{lemma_concern} and properties listed above, we have
\begin{equation}\label{eq_108}
 2| D_{\Omega,p}(\boldsymbol{U}^{2}(\boldsymbol{\Delta}^{2})^\top,\boldsymbol{\Delta}^{2}(\boldsymbol{U}^{2})^\top)| \leqslant 2\|\boldsymbol{\Omega}-p\boldsymbol{J}\| \sum_{i=1}^n \|\boldsymbol{U}_{i,\cdot}^{2}\|_2^2\|\boldsymbol{\Delta}_{i,\cdot}^{2}\|_2^2\leqslant 2\|\boldsymbol{\Omega}-p\boldsymbol{J}\|\nu_r  \|\boldsymbol{\Delta}\|_F^2
\end{equation}
and
\begin{equation}\label{eq_109}
\begin{split}
 4| D_{\Omega,p}(\boldsymbol{U}^{1}(\boldsymbol{\Delta}^{1})^\top,\boldsymbol{\Delta}^{2}(\boldsymbol{U}^{2})^\top)| \leqslant& 4\|\boldsymbol{\Omega}-p\boldsymbol{J}\|\sqrt{\sum_{i=1}^n \|\boldsymbol{U}_{i,\cdot}^{1}\|_2^2\|\boldsymbol{\Delta}_{i,\cdot}^{2}\|_2^2}\sqrt{\sum_{i=1}^n \|\boldsymbol{U}_{i,\cdot}^{2}\|_2^2\|\boldsymbol{\Delta}_{i,\cdot}^{1}\|_2^2}\\
 \leqslant & 2\|\boldsymbol{\Omega}-p\boldsymbol{J}\|\nu_r  \|\boldsymbol{\Delta}^{1}\|_F^2+ 2\|\boldsymbol{\Omega}-p\boldsymbol{J}\|\nu_r  \|\boldsymbol{\Delta}^{2}\|_F^2\\
 \leqslant & 2\|\boldsymbol{\Omega}-p\boldsymbol{J}\|\nu_r  \|\boldsymbol{\Delta}\|_F^2.
\end{split}
\end{equation}

Now putting estimations of terms in \eqref{eq_030} listed above together, i.e., \eqref{eq_105}, \eqref{eq_106}, \eqref{eq_107}, \eqref{eq_108} and \eqref{eq_109}, we have

\begin{equation}\label{eq_033}
\begin{split}
 &| D_{\Omega,p}(\boldsymbol{U}\boldsymbol{\Delta}^\top+\boldsymbol{\Delta}\boldsymbol{U}^\top,\boldsymbol{U}\boldsymbol{\Delta}^\top+\boldsymbol{\Delta}\boldsymbol{U}^\top) | \\
 \leqslant & \|\boldsymbol{\Omega}-p\boldsymbol{J}\|\left[ 3\nu_r  \|\boldsymbol{U}^{2}\|_F^2 +3\sum_{i=1}^n \|\boldsymbol{\Delta}_i\|_2^4 + 4\nu_r  \|\boldsymbol{\Delta}\|_F^2\right] + 10^{-4}p\|\boldsymbol{U}\boldsymbol{\Delta}^\top\|_F^2.
\end{split}
\end{equation}
Plugging estimations \eqref{eq_032}, \eqref{eq_061} and \eqref{eq_033} back to \eqref{eq_031}, we have  
\begin{equation*}
\begin{split}
K_2(\boldsymbol{X}) \leqslant &  | D_{\Omega,p}(\boldsymbol{\Delta}\boldsymbol{\Delta}^\top,\boldsymbol{\Delta}\boldsymbol{\Delta}^\top) | + 3 | D_{\Omega,p}(\boldsymbol{X}\boldsymbol{X}^\top-\boldsymbol{U}\boldsymbol{U}^\top,\boldsymbol{X}\boldsymbol{X}^\top-\boldsymbol{U}\boldsymbol{U}^\top)|  \\
\leqslant &  \|\boldsymbol{\Omega}-p\boldsymbol{J}\| \left[ 19\sum_{i=1}^n \|\boldsymbol{\Delta}_i\|_2^4 + 18\nu_r  \|\boldsymbol{\Delta}\|_F^2 + 9\nu_r  \sum_{i=s+1}^r \sigma_i \right] +3\times 10^{-4}p \|\boldsymbol{U}\boldsymbol{\Delta}^\top\|_F^2.
\end{split}
\end{equation*}

\end{proof}

\subsubsection{A proof of Lemma \ref{lemma_k4_psd}} 
\begin{proof}

 By the way we define $K_4(\boldsymbol{X})$ in \eqref{eq:K}, we have
 \begin{equation*}
	 \begin{split}
			K_4(\boldsymbol{X})  \leqslant &  | 6\langle  \boldsymbol{\Delta}\boldsymbol{\Delta}^\top ,\mathcal{P}_{\Omega}(\boldsymbol{N})\rangle - 6p \langle \boldsymbol{\Delta}\boldsymbol{\Delta}^\top,\boldsymbol{N}\rangle| + | 8\langle \boldsymbol{U}\boldsymbol{\Delta}^\top,\mathcal{P}_{\Omega}(\boldsymbol{N})\rangle - 8p\langle\boldsymbol{U}\boldsymbol{\Delta}^\top,\boldsymbol{N}\rangle|\\
			& + 6p\langle \boldsymbol{\Delta}\boldsymbol{\Delta}^\top ,\boldsymbol{N} \rangle\\
		 \leqslant &5\times 10^{-4}p \|\boldsymbol{\Delta}\boldsymbol{\Delta}^\top\|_F^2 +2\times 10^{-4}p \|\boldsymbol{U}\boldsymbol{\Delta}^\top\|_F^2  +C_2\frac{r\| \mathcal{P}_{\Omega} (\boldsymbol{N}) - p\boldsymbol{N} \|^2}{p}\\
		 & + 6p\langle \boldsymbol{\Delta}\boldsymbol{\Delta}^\top ,\boldsymbol{N} \rangle.
	 \end{split}
 \end{equation*}
 Here we use the fact that 
 \begin{equation*}
\begin{split}
 6|\langle \boldsymbol{\Delta}\boldsymbol{\Delta}^\top,\mathcal{P}_{\Omega} (\boldsymbol{N}) - p\boldsymbol{N}\rangle|\leqslant & 6 \frac{\sqrt{p}\|\boldsymbol{\Delta}\boldsymbol{\Delta}^\top\|_*}{\sqrt{r}}\frac{\sqrt{r}\|\mathcal{P}_{\Omega} (\boldsymbol{N}) - p\boldsymbol{N}\|}{\sqrt{p}}\\
 \leqslant& 6 \sqrt{p}\|\boldsymbol{\Delta}\boldsymbol{\Delta}^\top\|_F \frac{\sqrt{r}\|\mathcal{P}_{\Omega} (\boldsymbol{N}) - p\boldsymbol{N}\|}{\sqrt{p}}\\
 \leqslant & 5\times 10^{-4} p\|\boldsymbol{\Delta}\boldsymbol{\Delta}^\top\|_F^2 + C_2 \frac{r\| \mathcal{P}_{\Omega} (\boldsymbol{N}) - p\boldsymbol{N} \|^2}{p},
\end{split}
\end{equation*}		
where in the first line, we use matrix H\"{o}lder's inequality. For the second inequality, we use the fact that $\|\boldsymbol{\Delta}\boldsymbol{\Delta}^\top\|_* \leqslant \sqrt{r} \|\boldsymbol{\Delta}\boldsymbol{\Delta}^\top\|_F$. For the last inequality, we use the fact $2xy\leqslant wx^2 + \frac{y^2}{w}$ for all $w>0$. Use the same argument we also have

 \begin{equation*}
	 8|\langle \boldsymbol{U}\boldsymbol{\Delta}^\top, \mathcal{P}_{\Omega} (\boldsymbol{N}) - p\boldsymbol{N}\rangle| \leqslant      2\times 10^{-4} p\|\boldsymbol{U}\boldsymbol{\Delta}^\top\|_F^2 +C_2 \frac{r\| \mathcal{P}_{\Omega} (\boldsymbol{N}) - p\boldsymbol{N} \|^2}{p},
 \end{equation*}
which finishes the proof.

\end{proof}

\section{Discussions}
\label{sec:discussions}
This paper studies low-rank approximation of a positive semidefinite matrix from partial entries via nonconvex optimization. We established a model-free theory for local-minimum based low-rank approximation without any assumptions on its rank, condition number or eigenspace incoherence parameter. We have also improved the state-of-the-art sampling rate results for nonconvex matrix completion with no spurious local minima in \citet{ge2016matrix, ge2017no}, and have investigated the performance of the proposed nonconvex optimization in presence of large condition numbers, large incoherence parameters, or rank mismatching. The nonconvex optimization is further applied to the problem of memory-efficient Kernel PCA. Compared to the well-known Nystr\"{o}m methods, numerical experiments illustrate that the proposed nonconvex optimization approach yields more stable results in both low-rank approximation and clustering.

For future research, we are interested in understanding whether and how fast first-order methods converge to a neighborhood of the set of local minima with theoretical guarantees. In fact, a series of recent works in nonconvex optimization have discussed why and when first-order iterative algorithms can avoid strict saddle points almost surely. For example, in a very recent work by \citet{LPPSJ2017}, the authors show that under mild conditions of the nonconvex objective function, a variety of first order algorithms can avoid strict saddle points with almost all initialization, which extends the previous results in \citet{lee2016gradient} and \citet{PP17}. We are particularly interested in the robust version of the strict saddle points condition discussed in \citet{ ge2015escaping} and \citet{jin2017escape}, referred to as $(\theta, \gamma, \zeta)$-strict saddle, under which noisy stochastic/deterministic gradient descent methods are proven to converge to a neighborhood of the local minima. In fact, \citet[Theorem 12]{ge2017no} shows that the nonconvex optimization \prettyref{eq:obj_psd} satisfies certain $(\theta, \gamma, \zeta)$-strict saddle conditions as long as $\mtx{M}$ is exactly of rank $r$, its condition number and eigenspace incoherence parameter are well-bounded, and the sampling rate is sufficiently large, but their argument cannot be straightforwardly extended to the model-free settings. We plan to explore the $(\theta, \gamma, \zeta)$-strict saddle conditions for \prettyref{eq:obj_psd} under a model-free framework in future.

\section*{Acknowledgements}
We would like to acknowledge Taisong Jing for pointing us to the reference \citet{bhatia2013matrix}.

\bibliographystyle{plainnat}

    \bibliography{cite}

\begin{appendices}
	\section{Proof of Corollary \ref{coro_ex_psd_mu}}
	\begin{proof}
		The inequality \eqref{eq_100} gives $\|\boldsymbol{M}\|_{\ell_\infty} \leqslant \frac{\mu_r r \sigma_1}{n}$. Therefore, in the case $\rank(\boldsymbol{M})=r$, the approximation error bound \eqref{eq_062} becomes 
	\[ 
			 \left\|\widehat{\boldsymbol{X}}\widehat{\boldsymbol{X}}^\top - \boldsymbol{M}\right\|_F^2  \leqslant C_2    \sum_{i=1}^r \left\{\left[C_3  \left(\sqrt{\frac{n}{p}} +\frac{\log n}{p} \right)\frac{\mu_r r }{n}\sigma_1 -\sigma_i \right]_+\right\}^2.    
	\]
	Therefore, if 
	\[
			p\geqslant C \max\left\{\frac{\mu_r  r \kappa_r \log n}{n}, \frac{\mu_r^2  r^2\kappa_r^2 }{n}\right\}  
	\]
	with absolute constant $C$ sufficient large, we have 
	\[
		C_3  \left(\sqrt{\frac{n}{p}} +\frac{\log n}{p} \right)\frac{\mu_r r }{n}\sigma_1 \leqslant \sigma_i ,\quad i = 1,\cdots,r.
	\]
	In other words, $\widehat{\boldsymbol{X}}\widehat{\boldsymbol{X}}^\top = \boldsymbol{M}$. 
	
	Similarly, by definition \prettyref{eq:rong_spikiness}, in the case $\rank(\boldsymbol{M})=r$, we have 
		\begin{equation*}
			\|\boldsymbol{M}\|_{\ell_\infty} = \frac{\widetilde{\mu}_r^2 \operatorname{trace}(\boldsymbol{M})}{n}\leqslant \frac{\widetilde{\mu}_r^2 r \sigma_1}{n}.
		\end{equation*}
		Therefore, the approximation error bound \eqref{eq_062} becomes 
		\[ 
				 \left\|\widehat{\boldsymbol{X}}\widehat{\boldsymbol{X}}^\top - \boldsymbol{M}\right\|_F^2  \leqslant C_2    \sum_{i=1}^r \left\{\left[C_3  \left(\sqrt{\frac{n}{p}} +\frac{\log n}{p} \right)\frac{\widetilde{\mu}_r^2 r }{n}\sigma_1-\sigma_i \right]_+\right\}^2.
		\]
		Therefore, if 
		\[
			p\geqslant C \max\left\{ \frac{\widetilde{\mu}_r^2 r \kappa_r \log n}{n}, \frac{\widetilde{\mu}_r^4 r^2\kappa_r^2}{n} \right\} 
		\]
		with absolute constant $C$ sufficient large, we have $\widehat{\boldsymbol{X}}\widehat{\boldsymbol{X}}^\top = \boldsymbol{M}$.
	\end{proof}
	
	\section{Proof of Lemma \ref{lemma_k3_psd}}

	

	Here we present a proof of Lemma \ref{lemma_k3_psd}, this proof is exactly the proof in \citet{ge2017no} except keeping the extra negative term, we include the proof in \citet{ge2017no} here for completeness.
	\begin{proof}
	By \citet[Lemma 18]{ge2017no}, we have 
	\begin{equation}\label{eq_017}
		\begin{split}
			&\operatorname{vec}(\boldsymbol{\Delta})^\top \nabla^2 G_{\alpha}(\boldsymbol{X})\operatorname{vec}(\boldsymbol{\Delta}) - 4\langle \nabla G_{\alpha}(\boldsymbol{X}),\boldsymbol{\Delta}\rangle\\
			 =& 4\sum_{i=1}^n [(\|\boldsymbol{X}_{i,\cdot}\|_2-\alpha)_+]^3 \frac{\|\boldsymbol{X}_{i,\cdot}\|_2^2\|\boldsymbol{\Delta}_{i,\cdot}\|_2^2 -\langle \boldsymbol{X}_{i,\cdot},\boldsymbol{\Delta}_{i,\cdot}\rangle^2}{\|\boldsymbol{X}_{i,\cdot}\|_2^3}\\
			 &  + 12\sum_{i=1}^n [(\|\boldsymbol{X}_{i,\cdot}\|_2-\alpha)_+]^2 \frac{\langle \boldsymbol{X}_{i,\cdot},\boldsymbol{\Delta}_{i,\cdot}\rangle^2}{\|\boldsymbol{X}_{i,\cdot}\|_2^2} -16\sum_{i=1}^n [(\|\boldsymbol{X}_{i,\cdot}\|_2-\alpha)_+]^3 \frac{\langle \boldsymbol{X}_{i,\cdot},\boldsymbol{\Delta}_{i,\cdot}\rangle}{\|\boldsymbol{X}_{i,\cdot}\|_2}.
		\end{split}
	\end{equation}
	
		First of all, since we choose $\alpha \geqslant  100\sqrt{\nu_r} = 100 \|\boldsymbol{U}\|_{2,\infty}$, then for all $\boldsymbol{X}_{i,\cdot}$ satisfying $\|\boldsymbol{X}_{i,\cdot}\|_2\geqslant  \alpha$, we have
		\begin{equation}\label{eq_110}
			\langle\boldsymbol{X}_{i,\cdot},\boldsymbol{\Delta}_{i,\cdot} \rangle  = \langle\boldsymbol{X}_{i,\cdot}, \boldsymbol{X}_{i,\cdot}-\boldsymbol{U}_{i,\cdot}\rangle\geqslant  \|\boldsymbol{X}_{i,\cdot}\|_2^2 - \|\boldsymbol{X}_{i,\cdot}\|_2\|\boldsymbol{U}_{i,\cdot}\|_2 \geqslant  (1-0.01) \|\boldsymbol{X}_{i,\cdot}\|_2^2 \geqslant  0.99\|\boldsymbol{X}_{i,\cdot}\|_2^2,
		\end{equation}
		which gives an lower bound of the inner product between $\boldsymbol{X}_{i,\cdot}$ and $\boldsymbol{\Delta}_{i,\cdot}$, at the same time, we can also upper bound $\|\boldsymbol{\Delta}_{i,\cdot}\|_2$ by $\|\boldsymbol{X}_{i,\cdot}\|_2$:
		\begin{equation}\label{eq_111}
			\|\boldsymbol{\Delta}_{i,\cdot}\|_2 \leqslant \|\boldsymbol{X}_{i,\cdot}\|_2+\|\boldsymbol{U}_{i,\cdot}\|_2\leqslant 1.01\|\boldsymbol{X}_{i,\cdot}\|_2.
		\end{equation}	
	Plugging above two estimations \eqref{eq_110}, \eqref{eq_111} together with the fact that $|\langle \boldsymbol{X}_{i,\cdot},\boldsymbol{\Delta}_{i,\cdot}\rangle |^2\leqslant \|\boldsymbol{X}_{i,\cdot}\|_2^2\|\boldsymbol{\Delta}_{i,\cdot}\|_2^2$ into \eqref{eq_017}, we have
		\begin{equation}\label{eq_018}
			\begin{split}
				&\operatorname{vec}(\boldsymbol{\Delta})^\top \nabla^2 G_{\alpha}(\boldsymbol{X})\operatorname{vec}(\boldsymbol{\Delta}) - 4\langle \nabla G_{\alpha}(\boldsymbol{X}),\boldsymbol{\Delta}\rangle\\ 
				\leqslant & -15.68\sum_{i=1}^n [(\|\boldsymbol{X}_{i,\cdot}\|_2-\alpha)_+]^3\|\boldsymbol{X}_{i,\cdot}\|_2 +12\sum_{i=1}^n [(\|\boldsymbol{X}_{i,\cdot}\|_2-\alpha)_+]^2\|\boldsymbol{\Delta}_{i,\cdot}\|_2^2.\\
			\end{split}
		\end{equation}
	Moreover, for all $\boldsymbol{X}_{i,\cdot}$ satisfies $\|\boldsymbol{X}_{i,\cdot}\|_2\geqslant  5\alpha$, we can also upper bound $\|\boldsymbol{\Delta}_{i,\cdot}\|_2$ by $\|\boldsymbol{X}_{i,\cdot}\|_2$:
		\begin{equation}\label{eq_112}
			\|\boldsymbol{\Delta}_{i,\cdot}\|_2\leqslant \|\boldsymbol{X}_{i,\cdot}\|_2+\|\boldsymbol{U}_{i,\cdot}\|_2\leqslant 1.002 \|\boldsymbol{X}_{i,\cdot}\|_2,
		\end{equation}
		and also lower bound $\|\boldsymbol{X}_{i,\cdot}\|_2-\alpha$ by $\|\boldsymbol{\Delta}_{i,\cdot}\|_2$:
		\begin{equation}\label{eq_113}
			\|\boldsymbol{X}_{i,\cdot}\|_2-\alpha \geqslant  \left(1-\frac{1}{5}\right)\|\boldsymbol{X}_{i,\cdot}\|_2 \geqslant  \frac{400}{501} \|\boldsymbol{\Delta}_{i,\cdot}\|_2.
		\end{equation}
		Plugging \eqref{eq_112} and \eqref{eq_113} back to \eqref{eq_018}, we have
		\begin{equation*}
			\begin{split}
				&\operatorname{vec}(\boldsymbol{\Delta})^\top \nabla^2 G_{\alpha}(\boldsymbol{X})\operatorname{vec}(\boldsymbol{\Delta}) - 4\langle \nabla G_{\alpha}(\boldsymbol{X}),\boldsymbol{\Delta}\rangle\\
				\leqslant & 12\sum_{i,\|\boldsymbol{X}_{i,\cdot}\|_2< 5\alpha}[(\|\boldsymbol{X}_{i,\cdot}\|_2-\alpha)_+]^2\|\boldsymbol{\Delta}_{i,\cdot}\|_2^2 \\
				& + \left[ 12-15.68\times \frac{400}{501}\times \frac{1}{1.002}\right]\sum_{i,\|\boldsymbol{X}_{i,\cdot}\|_2\geqslant  5\alpha}[(\|\boldsymbol{X}_{i,\cdot}\|_2-\alpha)_+]^2\|\boldsymbol{\Delta}_{i,\cdot}\|_2^2\\
				\leqslant & 192\alpha^2 \|\boldsymbol{\Delta}\|_F^2 - 0.3\sum_{i,\|\boldsymbol{X}_{i,\cdot}\|_2\geqslant  5\alpha} \|\boldsymbol{\Delta}_{i,\cdot}\|_2^4\\
				\leqslant & 199.54\alpha^2 \|\boldsymbol{\Delta}\|_F^2 -0.3\sum_{i=1}^n \|\boldsymbol{\Delta}_{i,\cdot}\|_2^4,
			\end{split}
		\end{equation*}
		where the last inequality uses the fact that $\alpha \geqslant  100\sqrt{\nu_r}$.
	\end{proof}
	
\end{appendices}

\end{document}